\documentclass[10pt]{article}
\usepackage{latexsym}
\usepackage{fullpage}
\begin{document}
\input amssym.def
\input amssym
%
%
\font\tenbf=cmbx10  
\font\sevenbf=cmbx7
\font\fivebf=cmbx5
\newfam\bffam
\def\bold{\fam\bffam\tenbf}
\textfont\bffam=\tenbf
\scriptfont\bffam=\sevenbf
\scriptscriptfont\bffam=\fivebf
\font\tengoth=eufm10  
\font\sevengoth=eufm7
\font\fivegoth=eufm5
\newfam\gothfam
\def\goth{\fam\gothfam\tengoth}
\textfont\gothfam=\tengoth
\scriptfont\gothfam=\sevengoth
\scriptscriptfont\gothfam=\fivegoth
\font\bfgraec=cmmib10
%
%
\def\be{\begin{equation}}
\def\ee{\end{equation}}
\def\bea{\begin{eqnarray}}
\def\eea{\end{eqnarray}}
\def\bdm{\begin{displaymath}}
\def\edm{\end{displaymath}}
%
%
\def\a{\alpha}
\def\g{\gamma}
\def\G{\Gamma}
\def\d{\delta}
\def\D{\Delta}
\def\p{\varphi}
\def\t{\tau}
\def\s{\sigma}
\def\S{\Sigma}
\def\k{\kappa}
\def\l{\lambda}
\def\L{\Lambda}
\def\e{\beta}
\def\o{\omega}
\def\O{\Omega}
\def\z{\zeta}
\def\ra{\rangle}
\def\la{\langle}
\def\rb{\rbrack}
\def\lb{\lbrack}
\def\wt{\widetilde}
\def\A{\ifmmode \widehat{\bold A}
        \else $\widehat{\bold A}$ \fi}
\def\nb{\nabla}
\def\part#1{{\partial \over \partial x^{#1}}}
\def\qed{\ifmmode \quad\Box
         \else $\quad\Box$ \fi}
%
%
\def\K#1{\ifmmode {\Bbb K}^{#1} 
         \else $ {\Bbb K}^{#1}$ \fi}
\def\R{{\Bbb R}}
\def\C{{\Bbb C}}
\def\Q#1{\ifmmode {\Bbb Q}^{#1}\else $ {\Bbb Q}^{#1}$ \fi}
\def\Z#1{\ifmmode {\Bbb Z}^{#1}\else $ {\Bbb Z}^{#1}$ \fi}
\def\N{\ifmmode {\Bbb N}\else $ {\Bbb N}$ \fi}
\def\H{{\Bbb H}}
\def\RP#1{\ifmmode {\Bbb P}^{#1}(\R{})
          \else ${\Bbb P}^{#1}(\R{})$ \fi}
\def\CP#1{\ifmmode {\Bbb P}^{#1}(\C{})
          \else ${\Bbb P}^{#1}(\C{})$ \fi}
\def\HP#1{\ifmmode {\Bbb P}^{#1}(\H{})
          \else ${\Bbb P}^{#1}(\H{})$ \fi}
\def\KP#1{\ifmmode {\Bbb P}^{#1}(\K{})
          \else ${\Bbb P}^{#1}(\K{})$ \fi}
\def\spb{{\bold S}}
\def\spcb{\ifmmode \spb^c
     \else $ \spb^c $ \fi}
%
%
\def\SO{{\bold SO}}
\def\GL{{\bold GL}}
\def\O{{\bold O}}
\def\SU{{\bold SU}}
\def\U{{\bold U}}
\def\PU{{\bold PU}}
\def\Sp{{\bold Sp}}
\def\SL{{\bold SL}}
\def\spin{{\goth spin}}
\def\so{{\goth so}}
\def\su{{\goth su}}
\def\sp{{\goth sp}}
\def\sl{{\goth sl}}
\def\u{{\goth u}}
\def\Spin{{\bold Spin}}
\def\Spinc{\ifmmode {\bold Spin}^c
       \else ${\bold Spin}^c$ \fi}
\def\Spinh{\ifmmode {\bold Spin}^h
       \else ${\bold Spin}^h$ \fi}
\def\Cl{{\bold Cl}}
\def\CCl{{\Bbb C}{\bold l}}
\def\End{{\rm End}}
\def\Hom{{\rm Hom}}
\def\Aut{{\rm Aut}}
\def\Mat{{\rm Mat}}
\def\Bil{{\rm Bil}}
\def\Iso{{\rm Iso}}
\def\Sym{{\rm Sym}}
%
%
\def\zff{{\rm II}}
\def\im{{\rm im}}
\def\rg{{\rm rg}}
\def\id{{\rm id}\,}
\def\im{{\rm im}}
\def\ker{{\rm ker}}
\def\coker{{\rm coker}}
\def\det{{\rm det}}
\def\deg{{\rm deg}}
\def\sgn{{\rm sgn}}
\def\tr{{\rm tr}}
\def\dim{{\rm dim}}
\def\inf{\mathop{\rm inf}}
\def\ind{{\rm ind}}
\def\mod{{\rm mod}}
\def\grad{{\rm grad}}
\def\diag{{\rm diag}}
\def\trns{\phantom{}^t}
\def\lange{{\rm L\"ange}}
\def\span{{\rm span}}
\def\Ad{{\rm Ad}}
\def\ad{{\rm ad}}
\def\proj{{\rm proj}}
\def\max{{\rm max}}
\def\Ric{{\rm Ric}}
\def\rank{{\rm rank}\,}
\def\spa{{\rm span}}
\def\und{\qquad \mbox{and} \qquad}
%
%
\def\tupl#1#2{(#1_1,\dots,#1_{#2})}
\def\cst{{\rm const}}
\def\us{,\hskip-0.4mm,}
\def\os{``}
\def\buchst#1{\char\expandafter\number#1}
\def\addots{\mathinner{\mkern1mu\raise1pt\vbox{\kern7pt\hbox{.}}\mkern2mu
  \raise4pt\hbox{.}\mkern2mu\raise7pt\hbox{.}\mkern1mu}}
%
%
\def\kdirac{{\cal D}}
\def\kdiracq{\bar\kdirac}
\def\haupts{\hbox{\bfgraec\char27}}
\def\HB{{\bold H}}
\def\EB{{\bold E}}
%
%
\renewcommand{\theequation}{\thesection.\arabic{equation}}
\def\beweis{{\bf Beweis.}\hskip2mm}
\def\proof{{\bf Proof.}\hskip2mm}
\def\definition{{\bf Definition.}\hskip2mm}
\def\bemerkung{{\bf Bemerkung.}\hskip2mm}
\def\remark{{\bf Remark.}\hskip2mm}
\def\beispiel{{\bf Beispiel.}\hskip2mm}
\def\example{{\bf Example.}\hskip2mm}
\def\leer{\vskip 2.7mm \noindent}
\newtheorem{Satz}{Satz}[section]
\newtheorem{Proposition}{Proposition}[section]
\newtheorem{Theorem}{Theorem}[section]
\newtheorem{Corollary}{Corollary}[section]
\newtheorem{Korollar}{Korollar}[section]
\newtheorem{Lemma}{Lemma}[section]
\def\#{\sharp}
\def\b{\flat}
\def\Cr{{\rm Curv}}
\def\es{\lrcorner\;}
\def\ctr{\lrcorner}
\def\esn{\lrcorner_\circ}
\def\dn{\wedge_\circ}
\def\b{\flat}
\def\tfrac#1#2{{\textstyle\frac{#1}{#2}}}
\def\eqref#1{(\ref{#1})}
\title{Eigenvalue Estimates for the Dirac Operator on
 Quaternionic K\"ahler Manifolds}
\author{W. Kramer\thanks{Supported by the SFB 256 `Nichtlineare partielle
 Differentialgleichungen'},
 U. Semmelmann\thanks{Supported by the Graduiertenkolleg `Algebraische,
 analytische und geometrische Methoden und ihre Wechselwirkungen in der
 modernen Mathematik'}, 
  G. Weingart\thanks{Supported by the SFB 256 `Nichtlineare partielle
 Differentialgleichungen'}\\
 {\small Mathematisches Institut der Universit\"at Bonn}\\
 {\small Beringstra\ss{}e 1, 53115 Bonn, Germany}\\
 {\small\texttt{kramer@math.uni-bonn.de,
   uwe@math.uni-bonn.de, gw@math.uni-bonn.de}}}
\maketitle
\begin{abstract}
We consider the Dirac operator on compact quaternionic K\"ahler
manifolds and prove  a lower bound  for the spectrum.
This estimate is sharp since it is the first
eigenvalue of the Dirac operator on the quaternionic
projective space.
\end{abstract}

\tableofcontents

\section{Introduction}
On a compact Riemannian spin manifold $ ( M^n , \, g ) $
with positive scalar curvature $\k$
the eigenvalues of the Dirac operator $ D $ satisfy
$$
\lambda^2 \ge \frac{n}{n-1} \, \frac{\kappa_0}{4} \, ,
$$
where $\k_0$ is the minimum of the scalar curvature.
This estimate has been proven by Th.~Friedrich
(c.~f.~\cite{fri5}). This lower bound is sharp, since it is
attained as the first eigenvalue on the sphere.
A theorem of O.~Hijazi and A.~Lichnerowicz (\cite{hij2}, \cite{lic})
implies that the
case of equality cannot be attained on manifolds with non--trivial
parallel $k-$forms with $k\neq 0, n$.
There are two canonical classes of such manifolds
which in addition  have positive scalar curvature:
K\"ahler manifolds and quaternionic K\"ahler manifolds. The eigenvalue
estimate for manifolds of the first class has been improved
by K.--D.~Kirchberg in \cite{kir2} and \cite{kir6}
(see also \cite{lic1} and  \cite{hij1}). Again, this estimate is sharp:
the lower bound is equal to the first eigenvalue on complex projective
space in odd complex dimensions and on the product of $\CP{2m+1}$ with the
flat 2--torus in even complex dimensions.

The proof of the corresponding result
for the second class of manifolds is aim of the present article.
A quaternionic K\"ahler manifold is by definition an oriented
4n--dimensional Riemannian manifold with $n\geq 2$ 
whose holonomy group is contained in the subgroup 
$ \Sp(n)\cdot\Sp(1) \subset \SO(4n) $.
Equivalently they
are characterized by the existence of a certain parallel 4--form
$ \Omega $, the
so--called fundamental or Kraines form (c.~f.~\cite{bon}, \cite{krain}).
In even quaternionic dimensions $ n $ any quaternionic K\"ahler
manifold possesses a spin structure, whereas for odd $ n $ only
the quaternionic projective space is spin (c.~f.~\cite{sala}).
In this article we prove the following theorem:

\begin{Theorem}
Let $ ( M^{4n}, \, g ) $ be a compact quaternionic K\"ahler
spin manifold of positive scalar curvature $ \kappa $. Then
any eigenvalue $ \lambda $ of the Dirac operator satisfies
$$
\lambda^2 \ge \frac{n+3}{n+2} \, \frac{\kappa}{4} \,     .
$$
\end{Theorem}

We note that $ \kappa $  is indeed a constant, since any
quaternionic K\"ahler manifold is automatically Einstein.
This was first shown by {D.~V.~Alekseevskii} in \cite{alex2} and
\cite{alex2a} (see also \cite{ish3} and Lemma \ref{curvature}).

The spectrum of the Dirac operator $ D $ on the quaternionic
projective space has been computed by J.--L.~Milhorat in \cite{mil}.
The lower bound $ \tfrac{n+3}{n+2} \, \tfrac{\kappa}{4} $  turns out
to be the first eigenvalue  of $ D^2 $. Hence, the stated estimate is
sharp.

This lower bound was first conjectured by O.~Hijazi and
J.~--L.~Milhorat in \cite{hij3}. They gave first eigenvalue estimates
and proved the conjecture for quaternionic dimension $ n=2 $ and $ n = 3 $
(see also \cite{hij}).
Their approach was to define quaternionic K\"ahler twistor operators and to
use Weitzenb\"ock formulas to prove inequalities for the eigenvalues.
We will follow a similar approach by using representation
theory of the group $ \Sp(n)\cdot\Sp(1) $ to define natural twistor
operators.
The formalism we use systematically produces 
relations between differential operators of second order,
showing in particular that the Lichnerowicz--Weitzenb\"ock formula
alone is not sufficient to obtain the optimal estimate. 

In fact, we obtain a somewhat more general result. The spinor bundle of
a quaternionic K\"ahler manifold $ M^{4n} $ splits into a sum of
$ n $ subbundles and this decomposition is respected by the square of
the Dirac operator. We prove that the spectrum of $ D^2 $ on any
subbundle is bounded from below by the first eigenvalue of $ D^2 $
on the corresponding bundle on the quaternionic projective space.

The second author would like to thank D.~V.~Alekseevskii for many
interesting and helpful discussions. All of us would like to thank
C.~B\"ar and W.~Ballmann for encouragement and support.

\section{The Dirac Operator of a Quaternionic K\"ahler Manifold}
\subsection{Preliminaries on \textbf{Sp}$(n)$--Representations}

Let $ ( M^{4n}, \, g) $ be a quaternionic K\"ahler manifold.
The holonomy group $ \Sp(n) \cdot \Sp(1) \subset \SO(4n)$ reduces the
$ \SO(4n)$--bundle of orthonormal frames to a principal 
$ \Sp(n) \cdot \Sp(1) $--bundle $ P $, and the Levi--Civita connection on $M$
can be thought of to be given by a 1-form on $P$.
Any representation $ V $ of $ \Sp(n) \times \Sp(1) $ locally gives
a vector bundle ${\bold V}$ associated to $P$, which in addition does exist
globally provided the representation factors through $ \Sp(n) \cdot \Sp(1) $.

Let $H$ and ~$E$ be the defining complex representations of $\Sp(1)$ and 
$\Sp(n)$ with their invariant symplectic forms 
$ \sigma_H \in \Lambda^2 H^\ast$ and $\sigma_E \in \Lambda^2 E^\ast$ and
their compatible positive quaternionic structures $J$, e.~g.
$$
\begin{array}{rcl}
 J^2  & = & -1 \\
 \s_E(Je_1,Je_2) & = & \overline{\s_E(e_1,e_2)} \\
 \s_E(e,Je) & > & 0 \quad \hbox{for $e\neq 0$\, .}
\end{array}
$$
The symplectic form $\s_E$ defines an isomorphism $\#: E \to E^*,\;
e\mapsto e^\#:=\s_E(e,\,.\,)$ with inverse $ \flat : E^* \to E $.
Using either Gram's determinant or permanent the symplectic
form $  \sigma_E $   can be extended to    $ \Lambda^s E $
or $ \Sym^r E $. As an immediate consequence we obtain
$$
\sigma_E( e \wedge \eta_1, \, \eta_2 )
\; = \;
\sigma_E( \eta_1, \, e^{\sharp} \, \lrcorner \, \eta_2 ) \quad
\eta_1,\eta_2 \in  \Lambda^s E \,\,\hbox{resp.}\,\,\Sym^r E\, .
$$
The usual formula $\s_E(\,.\,,\,J.\,)$ then defines positive
definite hermitian forms, such that e.~g.~ $(e\wedge)^*=(Je)^\#\es$.
{\it Mutatis mutandis} analogous statements are true for $H$.

Let $ \{ e_i \} $  and $ \{ d e_i \} $ with $ d e_i ( e_j ) = \delta_{ij} $
be a dual pair of bases for $ E $, $E^\ast$ respectively. In terms of this
bases the symplectic form and its --- by $ \L^2 E \cong \L^2 E^* $
associated --- canonical bivector are given by
$$
\sigma_E
\; = \;
\tfrac{1}{2} \,  \sum d e_i \, \wedge \,   e_i^{\sharp}
\, \in \Lambda^2 E^*  \quad \qquad
L_E = \tfrac{1}{2} \,  \sum d e_i^\b \, \wedge \,   e_i \, 
\in \Lambda^2 E \, .
$$
Wedging with $ L_E $ determines a homomorphism
$ L  : \Lambda^{k-2} E \longrightarrow \Lambda^{k} E$
whereas contracting with $\s_E$ defines its adjoint
$ \L := L^\ast:\, \Lambda^{k} E \longrightarrow \Lambda^{k-2} E $.
The operators
$L$, $\L$ and $H:=\lb \L,L \rb$ fulfil the commutator rules of the
Lie algebra $\sl_2{\C}$ and, in addition, $H|_{\L^kE} = (n-k)\id$.
Therefore, $\L^kE= \im(L) \oplus \ker(\L)$ splits as 
$\Sp(n)$--representation and $\L_\circ^kE :=\ker(\L)$, the
{\it primitive} space, turns out to be irreducible. Hence,
the following decomposition is immediate:
$$
\Lambda^{q} E \quad = \quad
\bigoplus^{\lb\frac{q}{2}\rb}_{k=0}  \, \Lambda^{q-2k}_\circ  E, \qquad
0 \le q \le n \, .
$$
The primitive space is stable under
contraction with elements of $E^\ast$ but it is not preserved by the
wedge product. Therefore it is necessary 
to describe the projection 
$ e \dn \omega$ of
$ e \wedge \omega $ onto $ \; \Lambda^q_\circ  E $.
\begin{Lemma} \label{proj}
If $\omega \in   \Lambda^*_\circ E $ then
$e^\# \, \lrcorner \, \omega \; \in  \, \Lambda^*_\circ E$. Furthermore
$$
\begin{array}{rcl}
\Lambda(e \wedge \omega) & = &
e^{\sharp} \, \lrcorner \, \omega \\\\
e \wedge_\circ  \omega & = & e \wedge \omega  \;  - \;
   \tfrac{1}{n-k+1} \, L_E \wedge (e^{\sharp} \, \lrcorner \, \omega) \, .
\end{array}
$$
\end{Lemma}
To summarize properties of contraction and modified exterior
multiplication we have the following lemma:

\begin{Lemma} \label{kom2}
Let $\eta,\eta_i \in E^\ast$ and $e,e_i \in E$. Then
following relations are valid on $ \Lambda^s_\circ  E $:
$$\begin{array}{lclclcl}
 \{ \eta_1, \eta_2 \}  & = & 0 & \qquad \; &
 \{ \eta \, \lrcorner \, , e \, \dn  \,\} & = & \eta ( e )  +
  \tfrac{1}{n-s+1} \eta^{\flat} \, \dn  \,  e^{\sharp} \es \\\\
 \{ e_1  \dn  , e_2 \dn  \,\} & = & 0 & &
 \sum d e_i \, \lrcorner \, e_i \, \dn  &  = &
  (2n - s + 2 ) \tfrac{n-s}{n-s+1} \,  \id  \\\\
\sum e_i \, \dn  \, d e_i \, \lrcorner & = & s\,\id \, . &&&&
 \end{array}$$
\end{Lemma}
On $H$, there are similar equations which relate contraction and symmetric
product. It is convenient to modify the contraction on $ \Sym^r H $. 
For $\a \in H^\ast$, we define
$ \alpha \, \esn \, : \Sym^r H \rightarrow \Sym^{r-1} H $ by
$ \alpha \, \esn \, := \tfrac{1}{r} \, \alpha \, \lrcorner \,  $.
Let $ h \cdot $ denote the symmetric product with $ h \in H $.

\begin{Lemma}\label{kom1}
Let $h,h_i \in H$ and $\alpha,\alpha_i \in H^\ast$ Then the following
relations are valid on $ \Sym^r H $:
$$\begin{array}{lclclcl}
 \lb h_1 \cdot , h_2 \cdot \rb & = & 0 & \qquad\; &
 \lb\alpha \esn  , h \cdot \rb & = & -\frac{1}{r+1} \alpha^\b \cdot h^\# \esn 
  \\\\
 \lb \alpha_1 \esn  , \alpha_2 \esn  \rb & = & 0  & &
  \alpha(h)\id & = & h \cdot \alpha \esn  - \alpha^\b \cdot h^\# \esn  \\\\
 \sum h_i \cdot dh_i \esn   & = & \id & &
  \sum dh_i \esn  h_i \cdot & = & \frac{r+2}{r+1} \, \id \, .
\end{array}$$
\end{Lemma}

It follows from the Peter--Weyl theorem that any irreducible
$ \Sp(n) \times \Sp(1) $--module can be realized as a subspace
of $ H^{\otimes p} \otimes E^{\otimes q} $ for some $p$ and $q$. 
The representations
of $ \Sp(n) \cdot \Sp(1) $ are precisely those with $ p + q $ even.
Hence, any vector bundle on $ M $ associated to $ P $ can be expressed
in terms of the local bundles $ \HB $ and  $ \EB $.
All structures considered above carry over to the fibres of
the associated vector bundles.
For example, the complexified tangent bundle is defined by the representation
$ H \otimes E $, and we fix the identification
$$
T M^{\Bbb C } \; = \; \HB \otimes \EB       \, .
$$
Note that the real structure is simply $ \overline{h \otimes e}:= Jh \otimes
Je$. Finally, the Riemannian metric on $T M^{\Bbb C }$ is given by
$$
g_M \; = \; \sigma_H \otimes \sigma_E \, .
$$

\subsection{Spinor Bundle and Clifford Multiplication }

The aim of this section is to give an explicit formula for
the Clifford multiplication using the $H$--$E$--formalism.
The spinor module considered as  $\Sp(n) \times \Sp(1)$--representation
splits into a sum  of $n$ irreducible components. Hence, the spinor
bundle of a  $4n$--dimensional  quaternionic K\"ahler manifold
decomposes into a sum of $n$ subbundles which can be expressed
using the locally defined bundles $\EB$ and $\HB$. Likewise, this
decomposition can be defined by considering the Kraines form
$\Omega$ as endomorphism of the spinor bundle acting via
Clifford multiplication.

\begin{Proposition}
The spinor bundle $\spb(M)$ of a quaternionic K\"ahler manifold $M$
decomposes as
$$
\spb(M) \; = \; \bigoplus^n_{r=0} \, \spb_r(M)
  \; = \; \bigoplus^n_{r=0} \, \Sym^r \HB \otimes \Lambda^{n-r}_\circ \EB \, ,
$$
where each fibre of $\spb_r(M)$ is an eigenspace of  $ \; \Omega$  for the
eigenvalue
$$
\mu_r = 6n - 4r( r+2)  \, .
$$
The rank of the subbundle $\spb_r(M)$ is given by
$$
\rank (\spb_r(M)) = (r+1)\bigg( {2n \choose n-r } - 
 {2n \choose n-r -2}\bigg) \, .
$$
\end{Proposition}

The covariant derivative on $ \spb(M) $ induced by the Levi--Civita
connection on $ (M, g) $ respects the  decomposition
given above.
For further proceeding it is necessary to express the Clifford multiplication
on spinors in terms of the $H$--$E$--formalism.
\begin{Proposition}\label{cliff1}
For any tangent vector $ h \otimes e \in \HB \otimes \EB = TM^{\C}$,
the Clifford multiplication
$\mu( h \otimes e) : \spb(M) \rightarrow \spb(M)  $ is given by:
$$
 \mu(h \otimes e) \quad = \quad \sqrt{2} \;
 (
h \cdot \,  \otimes\, e^{\sharp} \, \lrcorner \,
\; + \;
h^{\sharp} \, \esn \, \otimes e \, \dn  \,)\, .
$$
In particular, the Clifford multiplication maps the subbundle
$\spb_r(M)$ to the sum $\spb_{r-1}(M) \oplus \spb_{r+1}(M)$.
\end{Proposition}
\leer
\proof
To check whether the above formula indeed defines the
Clifford multiplication it suffices to verify the general relation:
\be\label{cliff}
\mu(X) \circ \mu(Y) \; + \; \mu(Y) \circ \mu(X) \quad
 = \quad - \, 2 \; g(X, Y) \, .
\ee
We prove this for vectors $ X = h_1 \otimes e_1 $ and $ Y = h_2 \otimes e_2 $.
$$
\begin{array}{rcl}
\lefteqn{
\{
h_1^{\sharp} \, \esn \, \otimes \, e_1 \, \dn  \, , \,
h_2 \cdot \,  \otimes e_2^{\sharp} \, \lrcorner \,
\}
\;  = \;
[h_1^{\sharp} \, \esn \, , \, h_2 \cdot ]
\, \otimes \,
e_1 \, \dn  \,   e_2^{\sharp} \, \lrcorner \,
\, + \,
h_2 \cdot  h_1^{\sharp} \, \esn \,
\, \otimes \,
\{
e_1 \, \dn  \,  ,
e_2^{\sharp} \, \lrcorner \,
\}}&&         \\\\
& = &
- \,
\tfrac{1}{r+1} \,
h_1 \cdot  h_2^{\sharp} \, \esn \,
\, \otimes \,
e_1 \, \dn  \,  e_2^{\sharp} \, \lrcorner \,
\, + \,
h_2 \cdot  h_1^{\sharp} \, \esn \,
\, \otimes \,
(
\sigma_E(e_2, \, e_1)
\, + \,
\tfrac{1}{n-s+1} \,
e_2 \, \dn  \,   e_1^{\sharp} \, \lrcorner \,
)     \\\\
& = &
- \,
h_2 \cdot  h_1^{\sharp} \, \esn \,
\, \otimes \,
\sigma_E(e_1, \, e_2)
\, - \,
\tfrac{1}{r+1} \,
h_1 \cdot  h_2^{\sharp} \, \esn \,
\, \otimes \,
e_1 \, \dn  \,  e_2^{\sharp} \, \lrcorner \,
\, + \,
\tfrac{1}{n-s+1} \,
h_2 \cdot  h_1^{\sharp} \, \esn \,
e_2 \, \dn  \,   e_1^{\sharp} \, \lrcorner \,
\end{array}
$$
If we symmetrize over (1,2) and substitute $ s = n - r $
we obtain
$$
\begin{array}{rcl}
\{
\mu(h_1 \otimes e_1) , \,  \mu(h_2 \otimes e_2)
\}
& = &
- 2 \,
(
h_2 \cdot  h_1^{\sharp} \, \esn \,
\, - \,
h_1 \cdot  h_2^{\sharp} \, \esn \,
)
\, \otimes \,
\sigma_E(e_1, \, e_2)     \\\\
& = &
- 2 \, \sigma_H(h_1, \, h_2)  \, \sigma_E(e_1, \, e_2)  \,  \\\\
& = &
- 2 \, g( h_1 \otimes e_1 , \,  h_2 \otimes e_2 ) \,.
\end{array}
$$
The second assertion is clear from the definition. \qed
\leer
Thus, Clifford multiplication splits into two
components:
\be\label{mu}
\mu^+_- : \quad  TM \,\otimes\,\spb_r(M)\; \longrightarrow \; \spb_{r+1}(M)
\qquad \mbox{and} \qquad
\mu^-_+ : \quad  TM\,\otimes\,\spb_r(M)\; 
\longrightarrow \; \spb_{r-1}(M)  \, ,
\ee
where
$ \; \mu^+_-(e\,\otimes\, h\, \otimes\,\psi)
     = \sqrt{2} \, (h \cdot \otimes e \, \lrcorner \, ) \, \psi \; $
and
$ \; \mu^-_+ (e\,\otimes\, h\, \otimes\,\psi)
     = \sqrt{2} \,
     (h \, \esn \, \otimes e \, \dn  \,) \, \psi $.
We note that this definition makes sense also for $ \spb_r(M) $
replaced by $ \Sym^p \HB \otimes \Lambda^q_\circ \EB $. In this spirit
it is possible to define two operations similar to Clifford
multiplication:
$$
\begin{array}{rcl}
\mu^+_+: TM \otimes \Sym^p \HB \otimes \L_\circ^q \EB & \longrightarrow & 
\Sym^{p+1} \HB \otimes \L_\circ^{q+1} \EB \\ 
 h \otimes e \otimes \psi & \longmapsto &
{\sqrt 2}\, (h \cdot \, \otimes \,  e \, \dn\,)\psi \\\\
\hbox{and}&& \\\\
\mu^-_-: TM \otimes \Sym^p \HB \otimes \L_\circ^q \EB & \longrightarrow & 
\Sym^{p-1} \HB \otimes \L_\circ^{q-1} \EB \\
h \otimes e \otimes \psi &\longmapsto &
{\sqrt 2}\, (h \, \esn  \, \otimes \,  e \, \es \,)\psi \, .
\end{array}
$$

As a first application of Lemma \ref{cliff1} we will give a new
interpretation for the Kraines form acting on the subbundles $\spb_r(M)$
of the spinor bundle. For this we have to investigate the
Clifford multiplication with 2--forms. On an arbitrary spin manifold
it is defined by:
$$
\mu( X \wedge Y) =   \mu ( X) \mu (Y)  + g(X, \, Y)        \, .
$$
The space of 2--forms of a quaternionic K\"ahler manifold
split into the following components:
$$
\Lambda^2 TM \; \cong \; \Sym^2\HB \, \oplus \, \Sym^2 \EB  \, \oplus \,
(\Sym^2\HB \otimes \Lambda^2_\circ  \EB )  \, ,
$$
where $ \Sym^2\HB $ resp.~$ \Sym^2\EB $ is realized in $ \Lambda^2 TM $ as
$ \Sym^2\HB \otimes \sigma_E $ resp.~$\sigma_H \otimes \Sym^2\EB $.

\begin{Proposition}
Let $A$ be any element in $\Sym^2\HB \subset  \Lambda^2 TM $. Then
Clifford multiplication with $A$ on
$\spb_r(M) = \Sym^r\HB \otimes \Lambda^{n - r}_\circ \EB$, $( r \ge 1 ) $
is given by
$$
\mu (A) = 2 \, A \otimes \id \, ,
$$
where $ A $ acts as an element of $ \Sym^2H \cong \sp(1) $
on  $ \Sym^r\HB $. In particular, Clifford multiplication with
the Kraines form $ \Omega $ on $\spb_r(M)$ corresponds to the action of the
Casimir operator {\bold C} of $ {\goth sp}(1) $, i.~e.~ 
$$
\mu (\Omega) =  6n \, \id \; + \; 4 \, {\bold C} \,  \otimes \, \id
          \;   =  \; (6n \, - \,  4r(r+2) ) \id \, .
$$
\end{Proposition}

The next step is to introduce a scalar product
on the bundles $ \Sym^p\HB \otimes \Lambda^q_\circ  \EB $
which for $ q = n - p $ corresponds to the usual Hermitian
scalar product on the subbundle $ \spb_p(M) $ of the spinor bundle.
We have already seen that the symplectic forms $ \sigma_H $
and $ \sigma_E $  extend to positive definite Hermitean forms
on $ \Sym^r H $ and $ \Lambda_\circ ^s E $. Using this we define
on the tensor product $ \Sym^r \HB \otimes \Lambda_\circ ^s \EB $
the following twisted Hermitean form
$$
\langle
A_1 \otimes  \omega_1 , \,
A_2 \otimes  \omega_2
\rangle
\quad := \quad
\tfrac{1}{r!} \,
\sigma_H (  A_1, \,J A_2 ) \,   \sigma_E (  \omega_1 , \, J \omega_2 )   \, .
$$
The next lemma is an immediate consequence of this definition.

\begin{Lemma}              \label{skapro}
\begin{eqnarray}
\langle \mu^+_-(X) \psi_1, \, \psi_2\rangle
\quad & = & \quad
- \,  \langle \psi_1, \, \mu^-_+({\bar X}) \psi_2\rangle    \\[1ex]
\langle \mu^+_+(X) \psi_1, \, \psi_2\rangle
\quad & = & \quad
\langle \psi_1, \, \mu^-_-({\bar X}) \psi_2\rangle \, .
\end{eqnarray}
\end{Lemma}

\subsection{Dirac and Twistor Operators}
The aim of this section is to describe the decompositon of
$ TM \, \otimes \, \spb_r(M) $ and to introduce the corresponding
Dirac and twistor operators.
Using the Clebsch--Gordan formulas and similar formulas
for the fundamental representations of $\Sp(n)$
the tensor product decomposes for $ r \ge 1 $ as follows:

\bea 
TM \, \otimes \, \spb_r(M)
& \cong &
(\HB \otimes \EB) \, \otimes \, (\Sym^{r}\HB \otimes \Lambda^{n-r}_\circ  \EB) 
\nonumber \\[1.3ex]
& \cong &
  (\Sym^{r+1}\HB \otimes \Lambda^{n-r-1}_\circ  \EB)\;\oplus \; 
  (\Sym^{r-1}\HB \otimes \Lambda^{n-r+1}_\circ  \EB)
  \nonumber \\
& & \oplus \; (\Sym^{r+1}\HB \otimes \Lambda^{n-r+1}_\circ  \EB)
  \; \oplus \;
  (\Sym^{r-1}\HB \otimes \Lambda^{n-r-1}_\circ  \EB) \label{deco}\\
& & \oplus \; (\Sym^{r+1}\HB \otimes K^{n-r})
   \;  \oplus  \; (\Sym^{r-1}\HB \otimes K^{n-r})
  \nonumber \\[1.3ex]
& \cong & \spb_{r+1}(M) \;  \oplus \; \spb_{r-1}(M) \; \oplus \;
   ( \, S^+_{r} \; \oplus \; S^-_{r} \;\oplus \; V^+_r 
   \;\oplus\; V^-_r \,)\, .
   \nonumber
\eea
Here, $ K^{n-r} $ is the summand corresponding to the sum of the highest
weights in the decomposition of $ E \otimes \Lambda^{n-r}_\circ  E$.
Moreover, we introduced the notation
$S^{\pm}_{r} =  \Sym^{r \pm 1}\HB \otimes \Lambda^{n-r \pm 1}_\circ  \EB  \; $
and
$ \; V^{\pm}_r = \Sym^{r\pm 1}\HB \otimes K^{n-r} $.
In the case $ r = 0 \; $ and $r=n\;$ four of the
above summands vanish and we obtain:
\be\label{rzn}
(\EB \otimes  \HB) \, \otimes \, \Lambda^n_\circ \EB 
\; \cong \; \spb_1(M) \, \oplus \, V^+_0  \qquad \mathrm{and}\qquad
(\EB \otimes  \HB) \, \otimes \, \Sym^n\HB 
\; \cong \; \spb_{n-1}(M) \, \oplus \, S^+_n\, .
\ee

The two components of the Clifford multiplication define
natural projections onto the first two summands appearing in the
splitting (\ref{deco}). The remaining four summands constitute the kernel
of the Clifford multiplication. The projections onto $ S_r^+$
resp.~$S_r^-$ are given by $\mu^+_+$ resp.~$\mu^-_-$ defined in
the preceeding paragraph and, finally, the projectors onto $V^\pm$
will be denoted by $pr_{V^\pm}$. By applying these projectors on
the section $\nb \psi \in \G(TM \otimes \spb(M))$ we get the two
components of the Dirac operator
\be \label{defdirac}
D^{+}_{-} \,  := \,
\mu^{+}_{-} \circ \nabla : \spb_r(M)\,\longrightarrow\,\spb_{r +1}(M)\,\qquad 
D^{-}_{+} \,  := \,
\mu^{-}_{+} \circ \nabla : \spb_r(M) \, \longrightarrow\,\spb_{r-1}(M)\, ,
\ee
where $D = D^+_- + D^-_+$ is the full Dirac operator, and four twistor
operators:
\be
\begin{array}{rclrcl}
D^{+}_{+} &  := &
\mu^{+}_{+} \circ \nabla : \spb_r(M) \, \longrightarrow \,S_{r}^+ \, \qquad  &
D^{-}_{-} &  := &
\mu^{-}_{-} \circ \nabla : \spb_r(M) \, \longrightarrow \,S_{r}^- \,   \\\\
T^+ &  := &
pr_{V^+} \circ \nabla : \spb_r(M) \, \longrightarrow \, V^+  \qquad &
T^- &  := &
pr_{V^-} \circ \nabla : \spb_r(M) \, \longrightarrow \, V^- \,. 
\end{array}
\ee
The square of the Dirac operator respects the splitting of
the spinor bundle, i.~e.~$D^2 :  \spb_r(M) \, \longrightarrow \, \spb_r(M)$.
In particular, we have:
$ D^+_- D^+_- = 0 = D^-_+  D^-_+ $.
The adjoint operators are easily computed if one remembers the
scalar product introduced in the previous paragraph.

\begin{Lemma}  \label{dual}
$$
( D^+_- )^*  \;  =    \;    D_+^- \, ,  \quad \quad
( D^-_+ )^*  \;  =  \;      D_-^+ \, ,  \quad \quad
( D^+_+ )^*  \;  = \;   -   D_-^- \, ,   \quad \quad
( D^-_- )^*  \;  = \;  -    D_+^+ \, .
$$
\end{Lemma}
The proof is an easy consequence of Lemma \ref{skapro}.

\section{The Curvature Tensor on Quaternionic K\"ahler Manifolds}
\subsection{The Linear Space of Curvature Tensors}

The exterior algebra $\L V^*$ of a vector space $V$ defines a natural
multiplication $m:\Sym^2\L^2V^*\to\L^4V^*$. Its adjoint is usually
called comultiplication and is given by
$$
\Delta: \L^4V^\ast \to \Sym^2\L^2V^\ast,
\quad x\wedge y\wedge z\wedge w\mapsto
(x\wedge y)(z\wedge w)+(y\wedge z)(x\wedge w)+(z\wedge x)(y\wedge w) \, .
$$
The curvature tensor of a Riemannian manifold $M$ satisfies
Bianchi's first identity, i.e. it has to vanish on the image of $\Delta$.
Thus it is straightforward to call $\ker\,m$ the space of curvature tensors
on a vector space $V$. It is easier, however, to adopt the following
equivalent
\leer
\definition
The space of curvature tensors $\Cr\,V^*$ on $V$ is the subspace
$$\Cr\,V^*:=\spa\{ (\a\cdot\e)\times(\g\cdot\d):=(\a\wedge\g)(\e\wedge\d)+
 (\a\wedge\d)(\e\wedge\g)\,\,\mathrm{with}\,\,\a,\e,\g,\d\in V^*\}$$
of $\Sym^2\L^2V^*$. Note that the generators trivially satisfy
the `dual' Bianchi identity
$$(\a\cdot\e)\times(\g\cdot\d) + (\a\cdot\g)\times(\d\cdot\e) +
  (\a\cdot\d)\times(\e\cdot\g) = 0$$
and its variants using the apparent symmetries of
$(\a\cdot\e)\times(\g\cdot\d)$.

\begin{Lemma}
 There exists a natural isomorphism of vector spaces
 $$\begin{array}{ccc}
   \Sym^2\L^2V^*&\stackrel\cong\longrightarrow&\Cr\,V^*\oplus\L^4V^* \cr
    (\a\wedge\e)(\g\wedge\d) & \longmapsto & {1\over 3}
    \big((\a\cdot\g)\times(\e\cdot\d)-(\a\cdot\d)\times(\e\cdot\g)\big)
    \;\oplus\; \frac{1}{3}\a\wedge\e\wedge\g\wedge\d \, .
   \end{array}$$
 In the same vein
 $$\begin{array}{ccc}
   \Sym^2\Sym^2V^*&\stackrel\cong\longrightarrow &\Cr\,V^*\oplus\Sym^4\,V^*\cr
   (\a\cdot\e)(\g\cdot\d) & \longmapsto &
    {1\over 3}(\a\cdot\e)\times(\g\cdot\d)
    \;\oplus\; {1\over 3}\a\cdot\e\cdot\g\cdot\d \, .
   \end{array}$$
 This isomorphism reveals an alternative realisation of $\Cr\,V^*$ as the
 subspace of sectional curvature tensors in $\Sym^2\Sym^2V^*$.
\end{Lemma}
\leer
\proof
It is easy to write down an inverse of the first homomorphism using
the explicit formula
$$
(\a\wedge\e)(\g\wedge\d) =
  {1\over 3}((\a\cdot\g)\times(\e\cdot\d)-(\a\cdot\d)\times(\e\cdot\g))
  +\Delta({1\over 3}\a\wedge\e\wedge\g\wedge\d) \, .
$$
The inverse of the second homomorphism needs the analogues of $\Delta$
and $\times$:
$$
\begin{array}{rrcl}
 \Delta: & \Sym^4V &\longrightarrow & \Sym^2\Sym^2V \\
&  \a\cdot\e\cdot\g\cdot\d &\longmapsto & (\a\cdot\e)(\g\cdot\d)
  +(\a\cdot\g)(\d\cdot\e)+(\a\cdot\d)(\e\cdot\g)\\\\
 Cr^*: & \Sym^2\L^2V^* &\longrightarrow& \Sym^2\Sym^2V \\
&  (\a\wedge\e)(\g\wedge\d)
  & \longmapsto & (\a\cdot\g)(\e\cdot\d)-(\a\cdot\d)(\e\cdot\g) \, .
\end{array}
$$
With these definitions one finds
$$
(\a\cdot\e)(\g\cdot\d) = \Delta({1\over 3}\a\cdot\e\cdot\g\cdot\d)
  + Cr^*({1\over 3}(\a\cdot\e)\times(\g\cdot\d))\, . \qed
$$
\leer
\remark
Either of the above isomorphisms may be used to deduce
$$
\dim\,\Cr\,V^* = {N^2\over 12}(N^2-1) \qquad\rm{with}\,\,N=\dim\,V \, .
$$

In the next section the space $\Cr (H^*\otimes E^*)$ will be
decomposed into its irreducible $\SL\,H\times\SL\,E$-components.
For these calculations we need the additional

\begin{Lemma}
There exists a natural isomorphism
$$
  \begin{array}{ccc}
   \Sym^2V^*\otimes\L^2V^* & \stackrel\cong\longrightarrow &
   \L^2\Sym^2V^*\oplus\L^2\L^2V^* \cr \a\cdot\e\otimes\g\wedge\d
   & \longmapsto &
   {1\over 2}\big( (\a\cdot\g)\wedge(\e\cdot\d)+(\e\cdot\g)\wedge(\a\cdot\d)
     \;\oplus\; 
       (\a\wedge\g)\wedge(\e\wedge\d)+(\e\wedge\g)\wedge(\a\wedge\d)
    \big) \, .
   \end{array}
$$
\end{Lemma}
\leer
\proof
The partial inverse for the first summand is given by
$$
 (\a\cdot\g)\wedge(\e\cdot\d)\mapsto {1\over 2}\big(\a\cdot\e\otimes\g\wedge\d
   + \a\cdot\d\otimes\g\wedge\e + \g\cdot\e\otimes\a\wedge\d
   + \g\cdot\d\otimes\a\wedge\e \big) \, .
$$
Marginal changes are needed for the partial inverse of the second summand.
\qed
\leer
\begin{Lemma}\label{linjective}
Suppose that $V$ is a symplectic vector space with symplectic form $\s$. The
above isomorphism together with the embedding
$\Sym^2V^*\stackrel{\otimes\s}\longrightarrow\Sym^2V^*\otimes\L^2V^*$
induces mappings
$$
i_{\Sym}:\Sym^2V^*\to \L^2\Sym^2V^*\,,\qquad i_{\L}:\Sym^2V^*\to\L^2\L^2V^*\,,
$$
whereas $i_{\Sym}$ is always injective, $i_{\L}$ is
injective if and only if $\dim\,V\neq 2$.
\end{Lemma}
\leer
\proof
The symplectic form $\s=\sum_i de_i\otimes e_i^\#\in V^*\otimes V^*$
defines multiplication operators $\s\cdot$ and $\s\wedge$ on
$\Sym\,V^*\otimes\Sym\,V^*$ and $\L\,V^*\otimes\L\,V^*$ respectively,
which fit into commutative diagrams
\begin{center}
\begin{picture}(370,85)(0,-10)

\put(0,0){$V^\ast\otimes V^\ast$}
\put(0,60){$\Sym^2 V^\ast$}
\put(102,0){$\Sym^2V^\ast \otimes \Sym^2V^\ast$}
\put(120,60){$\L^2\Sym^2V^\ast$}
\put(42,4){\vector(1,0){55}}
\put(42,62){\vector(1,0){55}}
\put(17,55){\vector(0,-1){40}}
\put(142,55){\vector(0,-1){40}}
\put(63, 8){$\s\cdot$}
\put(59,68){$2i_{\Sym}$}

\put(220,0){$V^\ast\otimes V^\ast$}
\put(220,60){$\Sym^2 V^\ast$}
\put(322,0){$\L^2V^\ast \otimes \L^2V^\ast$}
\put(335,60){$\L^2\L^2V^\ast$}
\put(262,4){\vector(1,0){55}}
\put(262,62){\vector(1,0){55}}
\put(237,55){\vector(0,-1){40}}
\put(351,55){\vector(0,-1){40}}
\put(283,8){$\s\wedge$}
\put(279,68){$2i_{\L}$}
\end{picture}
\end{center}
because e.~g.~
\begin{eqnarray*}
  2i_{\Sym}(\a\cdot\e) & = & {1\over 2}\sum_i \lb (\a\cdot de_i)\wedge
  (\e\cdot e_i^\#) + (\e\cdot de_i)\wedge(\a\cdot e_i^\#)\rb \\
  & & \hookrightarrow\sum_i\lb de_i\cdot\a\otimes e_i^\#\cdot\e
  + de_i\cdot\e\otimes e_i^\#\cdot\a \rb\, .
\end{eqnarray*}
Associated to $\s$ is its canonical bivector $L$, which in turn defines 
contraction operators $L\es=\sum_i de_i^\b\es\otimes e_i\es$ in either case. 
The pairs
of operators $\s\cdot$ and $L\es$ or $\s\wedge$ and $L\es$ each generate
a $\sl_2$-algebra:
$$
\begin{array}{rcrcll}
 H_{\Sym} & := & \lb\s\cdot,-L\es\rb & = & \dim\,V+l+r &\mathrm{on}\quad
  \Sym^lV^*\otimes\Sym^rV^* \\ 
 H_{\L} & := & \lb L\es,\s\wedge\rb & = & \dim\,V-l-r &\mathrm{on}\quad
  \L^lV^*\otimes\L^rV^* \, . \\ 
\end{array}
$$
It remains to appeal to the representation theory of $\sl_2$ to conclude
that $\s\cdot$ is always injective, whereas
$\s\wedge: \L^lV^*\otimes\L^rV^*\to\L^{l+1}V^*\otimes\L^{r+1}V^*$
is injective only for $l+r<\dim\,V$.
\qed

\subsection{The Bianchi Identity for Tensor Products}

The aim of this section is to give a description of the space
$\Cr (H^*\otimes E^*)$ of curvature tensors on a tensor product
$H^\ast\otimes E^\ast$ in terms of the first Bianchi identity. To start with,
let's make explicit the usual isomorphisms
\begin{eqnarray*}
  \Sym^2(H^*\otimes E^*) & \stackrel\cong\longrightarrow &
    \Sym^2H^*\otimes\Sym^2E^* \oplus \L^2H^*\otimes \L^2E^*\\
  \L^2(H^*\otimes E^*) & \stackrel\cong\longrightarrow &
    \Sym^2H^*\otimes\L^2E^* \oplus \L^2H^*\otimes \Sym^2E^*
\end{eqnarray*}
setting e.~g.
$$
  (\a\otimes\e)\wedge(\g\otimes\d) \longmapsto {1\over 2}\lb
   \a\cdot\g\otimes\e\wedge\d\oplus\a\wedge\g\otimes\e\cdot\d \rb \, .
$$
Using all the isomorphisms above, the space $\Sym^2\L^2(H^*\otimes E^*)$
can be decomposed along the lines
\begin{eqnarray*}
\lefteqn{\Sym^2\L^2(H^*\otimes E^*)} && \\
 &\cong&\Sym^2(\Sym^2H^*\otimes\L^2E^*\oplus\L^2H^*\otimes\Sym^2E^*)\\
 &\cong&\Sym^2(\Sym^2H^*\otimes\L^2E^*) \;\oplus \; 
        (\Sym^2H^*\otimes\L^2E^*\otimes \L^2H^*\otimes\Sym^2E^*)\;\oplus\;
        \Sym(\L^2H^*\otimes \Sym^2E^*) \\
 &\cong&
 \begin{array}{rcl}
 (\Cr H^* &\otimes & \Cr E^*)_{\hbox to 0pt{{\tiny H}\hss}} \\
 \Cr H^* &\otimes & \L^4E^* \\
 \Sym^4 H^* &\otimes & \Cr E^* \\
 \Sym^4 H^* &\otimes & \L^4E^* \\
 (\L^2\Sym^2 H^* &\otimes & \L^2\L^2 E^*)_{\hbox to 0pt{{\tiny H}\hss}}
 \end{array}
 \oplus 
 \begin{array}{rcl}
 \L^2\Sym^2 H^* &\otimes & \L^2\Sym^2E^* \\
 (\L^2\L^2 H^* &\otimes & \L^2\Sym^2E^*)_{\hbox to 0pt{{\tiny M}\hss}} \\
 (\L^2\Sym^2 H^* &\otimes & \L^2\L^2E^*)_{\hbox to 0pt{{\tiny M}\hss}} \\
 \L^2\L^2 H^* &\otimes & \L^2\L^2E^* 
 \end{array}
 \oplus
 \begin{array}{rcl}
 (\Cr H^* &\otimes & \Cr E^*)_{\hbox to 0pt{{\tiny E}\hss}} \\
 \Cr H^* &\otimes & \Sym^4E^* \\
 \L^4 H^* &\otimes & \Cr E^* \\
 \L^4 H^* &\otimes & \Sym^4E^* \\
 (\L^2\L^2 H^* &\otimes & \L^2\Sym^2 E^*)_{\hbox to 0pt{{\tiny E}\hss}}
 \end{array}
\end{eqnarray*}
where summands occuring twice are labelled by their respective columns.
In exactly the same manner the space $\Sym^2\Sym^2 H^*$ may be decomposed.
With $\Cr (H^*\otimes E^*)$ being contained in either space a look at
differences and similarities leads to

\begin{Lemma}
In terms of the decomposition of $\Sym^2\L^2(H^*\otimes E^*)$
above the first Bianchi identity is equivalent to the five equations
\be
\begin{array}{lrcl}
{\rm I} \quad &
  pr_{(\Cr H^* \otimes \Cr E^*)_H}R &=&
  pr_{(\Cr H^* \otimes \Cr E^*)_E}R \\[1ex]
{\rm II} \quad &  pr_{\Sym^4H^* \otimes \L^4 E^*}R &=& 0 \\[1ex]
{\rm II'} \quad & pr_{\L^4H^* \otimes \Sym^4 E^*}R &=& 0 \\[1ex]
{\rm III} \quad &
  pr_{(\L^2\Sym^2 H^* \otimes \L^2\L^2 E^*)_H}R &=&
  pr_{(\L^2\Sym^2 H^* \otimes \L^2\L^2 E^*)_M}R \\[1ex]
{\rm III'} \quad &
  pr_{(\L^2\L^2 H^* \otimes \L^2\Sym^2 E^*)_E}R &=& 
  pr_{(\L^2\L^2 H^* \otimes \L^2\Sym^2 E^*)_M}R
\end{array}
\ee
\end{Lemma}
\leer
\proof
In the first place one has to check whether the space of solutions to
${\rm I}$--${\rm III}$ has the appropriate dimension. Setting
$N = \dim E^*$ and $M = \dim H^*$ a first step could be to calculate 
$$
  \dim(\Sym^4 E^* \oplus \Cr E^* \oplus \L^4 E^*) \,=\,
  {N+3 \choose 4} + {N^2\over 12}(N^2 -1) + {N \choose 4}
  \,=\, {N^2\over 6}(N^2 +5)\nonumber
$$
to find the dimension of the space of solutions to be
$$
\begin{array}{c}
 {\displaystyle\frac{N^2}{6}(N^2+5)\dim\,\Cr H^*+\frac{M^2}{6}(M^2+5)
  \dim\,\Cr E^*-\dim\,\Cr H^*\;\dim\,\Cr E^*} \\[1ex]
 {\displaystyle +\quad\frac{N^2}{4}(N^2-1)\frac{M^2}{4}(M^2-1)
  \qquad = \qquad \frac{N^2M^2}{12}(N^2M^2-1) \, .}
\end{array}
$$
Thus it is sufficient to show that the generators
$(\a \otimes \a')\cdot(\g \otimes \g')\times (\e \otimes \e')
\cdot(\d \otimes \d')$ of $\Cr (H^*\otimes E^*)$ in fact
satisfy equations ${\rm I}$--${\rm III}$, equations ${\rm II'}$
and ${\rm III'}$ are then inferred using the symmetry in $H$ and
$E$. To begin with, ${\rm II}$ is easy, because under the above
isomorphisms/projections the element
 $\big((\a \otimes \a')\wedge (\e \otimes \e')\big)
 \cdot\big((\g \otimes \g')\wedge (\d \otimes \d')\big)$
is mapped to
 ${1\over 72}\a\cdot\e\cdot\g\cdot\d \otimes
  \a'\wedge\e'\wedge\g'\wedge\d' \in\Sym^4 H^*\otimes \L^4E^*$.
Symmetrizing $ (\e \otimes \e') \leftrightarrow (\d \otimes \d')$
to get the image of the above generator yields 0. The argument for
equation ${\rm III.}$ is similar, yet slightly more complicated
as two images have to be calculated. Symmetrizing the two results
show that the generators are mapped to the same element in the
two copies of $\L^2\Sym^2H^*\otimes\L^2\L^2E^*$.

The trickiest part of the proof is ${\rm I}$, because the `dual' Bianchi
identity for the generators of $\Cr H^*$ and $\Cr E^*$ is used. Focussing
on the component $(\Cr\,H^*\otimes\Cr\,E^*)_H$, the chain of
isomorphism/projections maps
\begin{eqnarray*}
 \lefteqn{\bigl((\a \otimes \a')\wedge (\e \otimes \e')\bigr)
  \cdot\bigl((\g \otimes \g')\wedge (\d \otimes \d')\bigr)} & &\\
 &\mapsto& {1\over 4} (\a\cdot\e\otimes\a'\wedge\e')
   (\g\cdot\d\otimes\g'\wedge\d') \\
 &\mapsto& {1\over 8}(\a\cdot\e)(\g\cdot\d)\otimes
   (\a'\wedge\e')(\g'\wedge\d') \\
 &\mapsto& {1\over 72}(\a\cdot\e)\times(\g\cdot\d)\otimes\bigl(
   (\a'\cdot\g')\times(\e'\cdot\d')-(\a'\cdot\d')\times(\e'\cdot\g')\bigr)\,.
\end{eqnarray*}
The same element is projected to
${1\over 72}\bigl((\a\cdot\g)\times(\e\cdot\d)-(\a\cdot\d)\times(\e\cdot\g)
\bigr)\otimes(\a'\cdot\e')\times(\g'\cdot\d')$ in the other copy
$(\Cr\,H^*\otimes\Cr\,E^*)_E$. Even after symmetrization the results
look different, however the `dual' Bianchi identity for the generators
implies the following identity, which should fix this:
\begin{eqnarray*}
 \lefteqn{(\a\cdot\e)\times(\g\cdot\d)\otimes(\a'\cdot\g')\times(\e'\cdot\d')
  +(\a\cdot\d)\times(\g\cdot\e)\otimes(\a'\cdot\g')\times(\d'\cdot\e')} & &\\
 & = & - (\a\cdot\g)\times(\e\cdot\d)\otimes(\a'\cdot\g')\times(\e'\cdot\d')\\
 & = & (\a\cdot\g)\times(\e\cdot\d)\otimes(\a'\cdot\e')\times(\d'\cdot\g')
  +(\a\cdot\g)\times(\e\cdot\d)\otimes(\a'\cdot\d')\times(\g'\cdot\e') \, .
\qed
\end{eqnarray*}

\subsection{Solutions of the Bianchi Identity}

With the first Bianchi identity expressed in terms of the $H$--$E$--formalism
it is now easy to derive the decomposition of the curvature tensor of a
quaternionic K\"ahler manifold (compare \cite{alex2}, \cite{alex2a} or
\cite{sala}) into a linear combination of `curvature tensors' of the following
type:
\leer
\definition
Define the following $ \End(\HB\otimes\EB)$-valued 2-forms on $\HB\otimes\EB$:
\bea
  R^H_{h_1\otimes e_1,h_2\otimes e_2} &=&
  \s_E(e_1,e_2)(h_1h_2 \otimes \id_E) \nonumber\\
  R^E_{h_1\otimes e_1,h_2\otimes e_2} &=&
  \s_H(h_1,h_2)(\id_H \otimes e_1e_2) \nonumber\\
  R^{hyper}_{h_1\otimes e_1,h_2\otimes e_2} &=& \s_H(h_1,h_2)
  (\id_H \otimes {\goth R}_{e_1,e_2}),
\eea
where ${\goth R} \in \Sym^4E^*$ and ${\goth R}_{e_1,e_2}:
e\mapsto{\goth R}(e_1,e_2,e,.)^\b$ is an endomorphism of $E$.
\leer
\begin{Lemma}\label{curvature}
A quaternionic K\"ahler manifold $M$ is Einstein, and its
curvature tensor is given by
\be
  R = -\frac{\k}{8n(n+2)}(R^H + R^E) + R^{hyper} \, .
\ee
where $\k$ is the scalar curvature of $M$ and the symmetric
4-form ${\goth R}$ is necessarily the symmetrisation: 
$$
  {\goth R}(e_1,e_2,e_3,e_4) = \frac{1}{24 \s_H(h_1,h_2)\s_H(h_3,h_4)}
  \sum_{\t\in\,S_4} \la R_{{h_1}\otimes e_{\t 1},{h_2}\otimes e_{\t 2}}
   {h_3}\otimes e_{\t 3},{h_4}\otimes e_{\t 4}\ra \, .
$$
which is independent of the choice of the $h_i$ as long as 
$\s_H(h_1,h_2)\s_H(h_3,h_4) \neq 0$.
\end{Lemma}
\leer
\proof
By definition the curvature tensor is a 2-form on
$TM^\C\cong\HB\otimes\EB$ with values in $\sp H\oplus\sp E$,
i.~e.~it is an element of
$R \in \Sym^2(\Sym^2H^* \otimes \s_E \oplus \s_H\otimes \Sym^2E^*)$.
Of course $(\s_E)$ considered as an element of $\L(\L^2\EB)$ satisfies
$(\s_E)\wedge(\s_E)=0$, though $\s_E \wedge \s_E \neq 0$ for $\dim\, E\neq 2$.
Thus the $E$-symmetric part of $R$ can be written as
$$
  \frac{1}{2}(\s_H)^2 \otimes r^E \oplus \frac{1}{2}(\s_H)^2 \otimes {\goth R}
$$
for some $r^E\in\Cr E^*$ and ${\goth R} \in \Sym^4E^*$. Equation 
${\rm II'}$ is then trivially satisfied, because $\s_H\wedge\s_H=0$.
The $H$-symmetric part, however, reduces to 
$$
  r^H \otimes \frac{1}{2}(\s_E)^2
$$
for some $r^H\in\Cr H^*$, as a nonzero contribution from $\Sym^4H^*$
is incompatible with equation ${\rm II}$ with $\s_E\wedge\s_E$ being
nonzero.

Equation ${\rm III}$ then implies that the mixed part of $R$, 
an element of $\Sym^2H^*\otimes\s_E\otimes\s_H\otimes\Sym^2E^*$,
is mapped to $0$ under the projection to $\L^2\Sym^2H^*\otimes\L^2\L^2E^*$.
On the other hand, Lemma (\ref{linjective}) ensures this projection to be
injective on this particular subspace of
$\Sym^2H^*\otimes\L^2E^*\otimes\L^2H^*\otimes\Sym^2E^*$,
so that the mixed curvature part has to vanish.

Finally equation ${\rm I}$ shows that $r^E$ is a scalar multiple of
$pr_{\Cr\,E^*}({1\over 2}(\s_E)^2)$ and that $r^H$  is the same
multiple of $pr_{\Cr\,H^*}({1\over 2}(\s_H)^2)$, because this is
the only way to satisfy
$$
  r^E\otimes pr_{\Cr H^*}\Big(\frac{1}{2}(\s_H)^2 \Big) \, = \,
  pr_{\Cr E^*}\Big(\frac{1}{2}(\s_E)^2 \Big)\otimes r^H \, .
$$
Thus the first Bianchi identity implies that at any point of $M$
the curvature tensor of the quaternionic K\"ahler manifold is a
linear combination of $R^H+R^E$ and $R^{hyper}$ with:
\bea
  R^H &=& pr_{\Cr H^*}\Big( \frac{1}{2}(\s_H)^2 \Big)
  \otimes{1\over 2}(\s_E)^2 \nonumber\\
  R^E &=& \frac{1}{2}(\s_H)^2\otimes pr_{\Cr E^*}\Big(
  {1\over 2}(\s_E)^2 \Big) \nonumber\\
  R^{hyper} &=& \frac{1}{2}(\s_H)^2 \otimes {\goth R} \, .
\eea
To determine $R$ completely, it is convenient to calculate
the Ricci curvature of $M$, as its definition
$$
  \Ric(X,Y) = \tr(Z \mapsto R_{Z,X}Y), \quad X,Y,Z \in TM
$$
is easy to handle. Note that for real vectors $X$, $Y$ the endomorphism
$R_{.,X}Y$ is already defined over $\R$, so that its trace may be calculated
over $\R$ or $\C$. The contributions from the different components of $R$ are
\begin{itemize}
\item
$\Ric^{hyper}(h_1\otimes e_1,h_2\otimes e_2)$ is the trace of
the factorizable endomorphism
$$
  h\otimes e \mapsto \s_H(h,h_1)h_2 \otimes {\goth R}(e,e_1,e_2,.)^\b 
        = -(h^\#_1\otimes h_2)h \otimes {\goth R}_{e_1,e_2}e \, .
$$
Its trace is thus the product of the partial traces, however,
${\goth R}_{e_1,e_2} \in \Sym^2E^*\cong\sp E$ is trace-free.
\item
$\Ric^E(h_1\otimes e_1,h_2\otimes e_2)$ is the trace of the endomorphism
$$
  h\otimes e \mapsto \s_H(h,h_1)h_2 \otimes
  \big(\s_E(e,e_2)e_1 + \s_E(e_1,e_2)e\big)
  = -(h_1^\#\otimes h_2)h \otimes
  \big(-e_2^\#\otimes e_1 + \s_E(e_1,e_2)\id_E\big)e \, .
$$
which factorizes, too. Its trace is thus
$$
  \Ric^E(h_1\otimes e_1,h_2\otimes e_2) = -(2n+1)\s_H(h_1,h_2)\s_E(e_1,e_2)\, .
$$
\item
The same argument goes through for $\Ric^H$; the different dimensions
of $H$ and $E$ account for the slightly changed result:
$$
  \Ric^H(h_1\otimes e_1,h_2\otimes e_2) = -3\s_H(h_1,h_2)\s_E(e_1,e_2)\, .
$$
\end{itemize}
The Ricci curvature being a multiple of the metric the quaternionic
K\"ahler manifold $M$ is Einstein and the scalar curvature $\k$ is
thus constant on $M$. The coefficient of $R^H+R^E$ in $R$ is then
fixed by
\be
  \Ric(X,Y)= {\k\over 4n}g(X,Y) = -{\k\over 8n(n+2)}(\Ric^H+\Ric^E)(X,Y) \, .
\qed
\ee
\leer

At the end of this section a strange and remarkable feature of the
hyperk\"ahler part $ R^{hyper}$ of the curvature tensor should be pointed out.
This feature is convenient for deriving the Weitzenb\"ock formulas below,
nevertheless it implies that essentially new ideas are needed to discuss
the limit case in the eigenvalue estimate:
$R^{hyper}$ eludes all Weitzenb\"ock formulas!

To be more precise, the symmetric 4-form ${\goth R}$ induces a natural
endomorphism on every vector bundle associated to the principal
$\Sp(1)\Sp(n)$--bundle in the following way: 
\begin{eqnarray*}
  \Sym^4E^* & \stackrel\Delta\longrightarrow \,
  \Sym^2E^* \otimes \Sym^2E^*  \, \stackrel{m}\longrightarrow &
  {\cal U}(\sp E) \\ {\goth R} & \longmapsto & Q_{\goth R} \, .
\end{eqnarray*}
where ${\cal U}(\sp E)$ is the universal envelopping algebra of 
$\sp E$. This mapping is injective by Poincar\'e--Birkhoff--Witt,
nevertheless the endomorphism induced on the spinor bundle is
trivial:

\begin{Lemma}
Any symmetric 4-form in $\Sym^4E^*$ induces the trivial endomorphism on $\L E$.
\end{Lemma}
\leer
\proof
As the 4th powers span $\Sym^4E^*$, it is sufficient to check
the lemma only for an element of the form $\frac{1}{24}\a^4$. The
comultiplication maps it to ${1\over 2}\a^2\otimes{1\over 2}\a^2$.
Considered as an endomorphism ${1\over 2}\a^2$ denotes $A=\a\otimes\a^\b$,
extended as derivation to $\L E$. As such it satisfies $A^2=0$:
\begin{eqnarray*}
 A^2(e_1\wedge\ldots\wedge e_s) & = &
  \sum_{i=1}^s e_1\wedge\ldots\wedge A^2e_i \wedge\ldots\wedge e_s\\
 & & \qquad +2\sum_{i<j} e_1\wedge\ldots\wedge Ae_i\wedge\ldots\wedge
     Ae_j\wedge \ldots\wedge e_s \, .
\end{eqnarray*}

In fact, the first sum is zero, because already $A^2=0$ on $E$, and the
second vanishes identically, because the endomorphism $A$ of $E$ has rank 1.
\qed
\leer

\begin{Corollary}
The comultiplication of ${\goth R}$ may be written as
$$\Delta {\goth R}=
 {1\over 2}\sum_{i,j}de_i\cdot de_j\otimes{\goth R}(e_i,e_j,.,.)$$
such that
\be
 \label{qzero}
 \sum_{i,j}\id\otimes(de_j^\b\dn de_i\es +de_i^\b\dn de_j\es )
  {\goth R}_{e_i,e_j} = 0
\ee
holds as an operator identity on the spinor bundle.
\end{Corollary}

\section{The Universal Weitzenb\"ock Formula}

In this section we develop the main tool for  
the proof of the lower bound. It turns out that all necessary information
can be encoded in a single matrix. To motivate the further proceeding,
it is useful to remember the situation in the Riemannian case. For 
deriving the eigenvalue estimate in that case, two things are needed: the
Lichnerowicz--Weitzenb\"ock formula and the twistor operator.

\subsection{The Riemannian Estimate Revisited}

Let $(M^n,g)$ be Riemennian spin manifold with tangent resp.~spinor bundle
$TM $ resp.~$\spb(M)$ associated to the representations $V$ and $\S$ of
$\Spin(V)$. The associativity of the tensor product 
$(V \otimes V) \otimes \S \cong V \otimes (V \otimes \S) $ can be thought of
as two different ways of decomposing $ V \otimes V \otimes \S $ into
irreducible $\Spin(V)$--representations; the isomorphism induced is
expressed by an invertible matrix ${\cal W}$ for each isotypical component.

Accordingly, the special section 
$\nb^2\psi \in \G(TM\otimes TM \otimes \spb(M))$ splits into sections of the
associated bundles inducing two different sets of 2nd order differential
operators which are related by the same matrix.
We will only be interested in differential operators from the spinor
bundle to itself.

In the Riemannian case $V\otimes V\otimes \S$ contains two copies of $\S$, such
that the associativity as above is expressed by a $2\times 2$--matrix. Let's
look at $(V \otimes V) \otimes \S$ first. Due to the obvious decomposition
$$
  V \otimes V \cong \L^2V \oplus \Sym^2_\circ V \oplus \C \, ,
$$
the projection
onto the two summands in $(V \otimes V) \otimes \S$ is given by contracting of
the two $V$ factors with the metric $g$ and the operation $ \L^2 V \cong 
\spin(V)$ on $\S$. The space $\Sym^2_\circ V \otimes \S$ contains no copy of
$\S$. To summarize:
$$
\begin{array}{rrcl}
  pr_{\C}: & V \otimes V \otimes \S & \longrightarrow & \S \\
&  e_1\otimes e_2\otimes \psi  & \longmapsto & g(e_1,e_2)\psi \, ,\\\\
  pr_{\L^2}: &  V \otimes V \otimes \S &\longrightarrow& \S \\
  & e_1\otimes e_2\otimes \psi & \longmapsto& e_1e_2\psi + g(e_1,e_2)\psi \, .
\end{array}
$$
Strictly speaking, there is a factor $\frac{1}{2}$ missing in $pr_{\L^2}$, but
this is of no concern at this point.
On the other hand, the projection to the first $\S$ in $V \otimes (V\otimes\S)$
is merely twice the Clifford multiplication. The complementary projection
involves the kernel $K \subset V\otimes \S$ of the Clifford multiplication
and is contraction of the two $V$--factors in $V \otimes K$:
$$
\begin{array}{rrclcl}
  pr_\S: & V \otimes V \otimes \S & \longrightarrow & \S &&\\
  & e_1 \otimes e_2 \otimes \psi & \longmapsto & e_1e_2\psi \, , && \\\\
  pr_K: & V \otimes V \otimes \S & \longrightarrow& 
    V\otimes K \subset V \otimes V \otimes \S & \longrightarrow & \S \\
  & e_1 \otimes e_2 \otimes\psi & \longmapsto & e_1 \otimes (e_2\otimes
    \psi + \frac{1}{n} \sum_i e_i\otimes e_i\cdot e_2\cdot\psi) &
    \longmapsto & g(e_1,e_2)\psi + \frac{1}{n} e_1e_2\psi \, .
\end{array}
$$

With these calculations it is now easy to provide a first example of a
matrix ${\cal W}$ as above, called the Weitzenb\"ock matrix of a Riemannian
manifold. It is the matrix appearing in the explicit formula for the
isomorphism induced by associativity:
$$
  \left(\begin{array}{l}pr_{\C} \\\\ pr_{\L^2}\end{array}\right) = 
  \left(\begin{array}{ll}-\frac{1}{n} & 1 \\\\
                    \frac{n-1}{n} & 1 \end{array}\right)
  \left(\begin{array}{l}pr_\S \\\\ pr_K\end{array}\right) \, .
$$

The last thing to do is to identify the projectors with differential operators
like $\nb^\ast \nb$ and $ D^2$ on the manifold. For a section
$\psi \in \G(\spb(M))$ the section $\nb^2\psi$ can be expanded in a sum over an
orthonormal basis leading to the following images of the projections
onto $\S$: 
$$
\begin{array}{rcl}
  pr_\C(\nb^2\psi) &=& \sum_{i,j} g(e_i,e_j) \nb^2_{e_i,e_j}\psi 
     = -\nb^\ast\nb \psi\\\\
  pr_{\L^2}(\nb^2\psi) &=& \sum_{i,j} (e_ie_j + g(e_i,e_j))\nb^2_{e_i,e_j}\psi 
     = \frac{\k}{4} \psi
\end{array}
$$
and 
$$
\begin{array}{rcl}
  pr_\S(\nb^2\psi) &=& \sum_{i,j} e_ie_j\nb^2_{e_i,e_j}\psi 
     = D^2 \psi \\\\
  pr_K(\nb^2\psi) &=& \sum_{i,j}\big(g(e_i,e_j)
    +\frac{1}{n}e_ie_j\nb^2_{e_i,e_j}\big)\psi = -T^\ast T \psi \, .
\end{array}
$$
The only projection which does not lead to a differential operator is the
well--known curvature term appearing in the Lichnerowicz--Weitzenb\"ock
formula. The final form of the Riemannian version of the Weitzenb\"ock
matrix formula reads:
\be \label{rwformel}
  \left(\begin{array}{c} -\nb^\ast\nb\psi \\\\ \frac{\k}{4}\psi \end{array}
  \right) = \left(\begin{array}{cc}-\frac{1}{n} & 1 \\\\
                    \frac{n-1}{n} & 1 \end{array}\right)
  \left(\begin{array}{c} D^2\psi \\\\ -T^\ast T\psi \end{array}\right) \, .
\ee

This matrix equation is a nice tool to produce some well--known formulas by
multiplying it with row vectors from the left, e.~g.~multiplying with
$(-1 \;\; 1)$ yields the Lichnerowicz--Weitzenb\"ock formula.
Taking the $L^2$-product with $\psi$, we get its integrated version:
\be 
  \left(\begin{array}{c} -\|\nb\psi\|^2 \\\\ \frac{1}{4}\la\psi,\k\psi\ra 
        \end{array}\right) = 
  \left(\begin{array}{cc}-\frac{1}{n} & 1 \\\\
                    \frac{n-1}{n} & 1 \end{array}\right)
  \left(\begin{array}{c} \|D\psi\|^2 \\\\  
   -\| T\psi\|^2 \end{array}\right) \, .
\ee
Estimating the twistor operator term to zero, the second row immediately gives
the well-known estimate of Friedrich \cite{fri5} for the first eigenvalue
of the Dirac operator. This rather simple example of the Riemannian
case shows the principle by which the eigenvalue estimate in
the quaternionic case will be obtained.

\subsection{The Weitzenb\"ock Matrix}

The same strategy is pursued in the quaternionic K\"ahler case, but under
holonomy $\Sp(1)\cdot\Sp(n)$ the splitting of $TM\otimes TM\otimes\spb(M)$
is much more delicate. Like in the Riemannian case, the splitting of
$ \nb^2\psi$, considered as section
in $TM \otimes TM \otimes \spb(M)$ yields two different sets of differential
operators related by a matrix ${\cal W}$. To calculate it, we have to look
at the splitting of
$$
  V \otimes V \otimes \S = 
  (H \otimes E) \otimes (H \otimes E) \otimes \bigoplus_{r=0}^n \big(\Sym^rH
\otimes
  \L_\circ^{n-r}E\big)
$$
into irreducibles, find the $\S$-summands and determine how they are 
embedded in $V \otimes V \otimes \S$.

There is a simple observation which makes the following calculations feasible.
All representations occuring are products of representations of $\Sp(1)$ or
$\Sp(n)$. Hence it suffices to look at the $H$--part resp.~the $E$--part of
the splitting, that means, one can look at the spaces 
$H \otimes H \otimes \Sym^rH$ and $E \otimes E \otimes \L_\circ^{n-r}E$ 
separately to
search for their $\Sym^rH$-- resp.~$\L_\circ^{n-r}E$--summands.

First, let's consider the $H$--part. The two sides of the isomorphism
$$
  (H \otimes H) \otimes \Sym^rH \cong H \otimes (H \otimes \Sym^rH)
$$
need to be decomposed, where on the left side $H \otimes H$ is acting as
endomorphisms on $\Sym^rH$ and on the right side $H$ is acting by its part
of the Clifford multiplication. On the right side the situation is easy:
$$
\begin{array}{rrclcl}
  pr_{-+}: & H \otimes H \otimes \Sym^rH & \longrightarrow &
  H \otimes \Sym^{r+1}H &\longrightarrow &  \Sym^rH \\\
  & h_1 \otimes h_2 \otimes s & \longmapsto & h_1 \otimes h_2 \cdot s
    & \longmapsto & h_1^\# \esn(h_2 \cdot s) \, , \\\\
  pr_{+-}: & H \otimes H \otimes \Sym^rH & \longrightarrow &
  H \otimes \Sym^{r-1}H & \longrightarrow & \Sym^rH \\\
  & h_1 \otimes h_2 \otimes s & \longmapsto & h_1 \otimes h_2^\# \esn s
   & \longmapsto & h_1 \cdot(h_2^\# \esn s) \, .
\end{array}
$$
Using Lemma \ref{kom1} the expression for the first projection
can be rewritten:
$$
  pr_{-+}(h_1 \otimes h_2 \otimes s) = 
  \s_H(h_1,h_2)s+ \frac{r}{r+1} h_1\cdot h_2^\#\esn s \, .
$$
On the left hand side, we use the isomorphism
$$
  H \otimes H \cong \Sym^2H \oplus \L^2H \cong \Sym^2H \oplus \C \, ,
$$
where $\C$ denotes the trivial representation, spanned by $\s_H$.
Hence,
$$
\begin{array}{rcl}
  pr_{\C}: H \otimes H \otimes \Sym^rH & \longrightarrow & \Sym^rH \\
  h_1 \otimes h_2 \otimes s & \longmapsto & \s_H(h_1,h_2) s \, .
\end{array}
$$
Elements of $\Sym^2H$ operate in the obvious manner as endomorphims
on $H$, namely as $ (h_1h_2)(h) = \s_H(h_1,h)h_2 + \s_H(h_2,h)h_1$,
extended as derivation to the symmetric algebra:
$$
\begin{array}{rcl}
  pr_{\Sym^2H}: H \otimes H \otimes \Sym^rH & \longrightarrow & \Sym^rH 
\nonumber\\
  h_1 \otimes h_2 \otimes s & \longmapsto & (h_1h_2)(s) \, .
\end{array}
$$
This can be written explicitely:
$$
  (h_1h_2)(s) = (h_2\cdot h_1^\#\es\; + h_1\cdot h_2^\#\es \;) s
    = r \s_H(h_1,h_2) s + 2h_1\cdot h_2^\#\es s = r \s_H(h_1,h_2) s + 2r 
h_1\cdot h_2^\#\esn s \, .
$$
To get an overview of this situation, a diagram will be helpful. The notation 
$\lb.\rb_{\Sym^rH}$ refers to the irreducible $\Sym^rH$-summand of the
space in question.

\begin{center}
\begin{picture}(400,150)(-155,-40)

\put(-140,100){$\lb\C\otimes\Sym^rH\rb_{\Sym^rH}$}
\put(-140,85){$\s_H(h_1,h_2) s$}
\put(135,100){$\Sym^rH$}
\put(135,85){$\s_H(h_1,h_2)s+ \frac{r}{r+1} h_1\cdot h_2^\#\esn s$}
\put(0,40){$H \otimes H \otimes \Sym^rH$}
\put(0,25){$h_1 \otimes h_2 \otimes s$}
\put(135,-20){$\Sym^rH$}
\put(135,-35){$h_1\cdot h_2^\#\esn s$}
\put(-155,-20){$\lb\Sym^2H\otimes\Sym^rH\rb_{\Sym^rH}$}
\put(-155,-35){$r \s_H(h_1,h_2) s + 2r h_1\cdot h_2^\#\esn s$}
\put(-10,45){\vector(-1,1){40}}
\put(-10,30){\vector(-1,-1){40}}
\put(80,45){\vector(1,1){40}}
\put(80,30){\vector(1,-1){40}}
\put(-45,60){$pr_{\C}$}
\put(-70,10){$pr_{\Sym^2H}$}
\put(105,60){$pr_{-+}$}
\put(105,10){$pr_{+-}$}

\end{picture}
\end{center}
From this diagram the $H$--part of the Weitzenb\"ock matrix is evident:
\be
  {\cal W}_H = \left(\begin{array}{cc}
   1 & -\frac{r}{r+1} \\\\ r & \frac{r(r+2)}{r+1} \end{array}\right) \, .
\ee

Turning to the $E$--part of the Weitzenb\"ock matrix, one has to look
at the two sides of the obvious isomorphism again:
\be \label{gl2}
  (E \otimes E) \otimes \L_\circ^{n-r}E 
\cong E \otimes (E \otimes \L_\circ^{n-r}E) \, .
\ee
In difference to the $H$--part there exists a kernel $K{n-r}$ of 
multiplication and contraction on the right side:
$$
  E \otimes \L_\circ^{n-r}E \cong \L_\circ^{n-r+1}E 
\oplus \L_\circ^{n-r-1}E \oplus K^{n-r} \, .
$$
Therefore, in the decomposition of $E \otimes (E \otimes \L_\circ^{n-r}E)$
the representation
$\L_\circ^{n-r}E$ occurs three times, namely twice by the $E$--part of the 
Clifford multiplication and in addition as an irreducible summand of
$E \otimes K^{n-r}$. The following lemma gives an explicit expression for
the projection onto $K^{n-r}$:
\begin{Lemma}
$$
\begin{array}{rrcl}
  \widetilde{pr}_{K}: & E \otimes \L_\circ^{n-r}E & \longrightarrow & K^{n-r}\\
  & e \otimes \omega & \longmapsto & e \otimes \omega 
    -\frac{1}{n-r+1} \sum_i e_i \otimes de_i \es e \dn \omega \\\\
   &&&  -\frac{r+2}{(n+r+3)(r+1)}\sum_i de_i^\b \otimes e_i \dn e^\# \es \omega
   \, .
\end{array}
$$
\end{Lemma}
\leer
\proof
Simple use of the number operators introduced in Lemma \ref{kom2} immediately
shows that the projection of $e \otimes \omega$ onto the kernel
of the $\dn$-product is $e \otimes \omega -\frac{1}{n-r+1} \sum_i e_i 
\otimes de_i \es e \dn \omega $. By the same argument $e\otimes \omega$
plus the third summand is the projection onto the kernel of contraction:
$$
\begin{array}{rcl}
  \sum_i de_i \es e_i \dn e^\# \es \omega 
&=& \sum_i \big( \s_E(de_i^\b, e_i) - e_i \dn de_i \es + \tfrac{1}{r+2}
   de_i^\b \dn e_i^\# \es \big) e^\# \es \omega \\\\
&=& ( 2n -(n-r-1) - \tfrac{n-r-1}{r+2})e^\# \es\omega \\\\
&=& \tfrac{(n+r+3)(r+1)}{r+2}e^\# \es\omega \, .
\end{array}
$$
It remains to show that the second resp.~the third summand already lies in
the kernel of contraction resp.~multiplication. This is clear for the second
summand because $\s_E \es e^\# \dn \omega = 0$ by definition. Upon
$\dn$--multiplication, the third summand reduces to the projection of
$L_E \wedge e^\# \es \omega$ onto $\L_\circ E$ which is zero by definition,
too.
\qed
\leer
Tensoring once again, the projection onto the irreducible 
$\L_\circ^{n-r}E$--summand
contained in $ E \otimes K^{n-r} \subset E \otimes E 
\otimes \L_\circ^{n-r}E$ can be
thought of as contracting the two $E$--factors.
All ingrediences for computing the right hand side are now at hand:
$$
\begin{array}{rrclcl}
  pr_{-+}: & E \otimes E \otimes \L_\circ^{n-r}E & \longrightarrow & 
  E \otimes \L_\circ^{n-r+1}E & \longrightarrow & \L_\circ^{n-r}E \\
  & e_1 \otimes e_2 \otimes \omega & \longmapsto &
  e_1 \otimes e_2 \dn \omega & \longmapsto & e_1^\# \es e_2 \dn \omega\, , 
\\\\
  pr_{+-}: & E \otimes E \otimes \L_\circ^{n-r}E & \longrightarrow & 
  E \otimes \L_\circ^{n-r-1}E & \longrightarrow & \L_\circ^{n-r}E \\
  & e_1 \otimes e_2 \otimes \omega & \longmapsto &
  e_1 \otimes e_2^\# \es \omega & \longmapsto & e_1 \dn e_2^\# \es \omega\, , 
\\\\
  pr_K: & E \otimes E \otimes \L_\circ^{n-r}E & 
\longrightarrow & E \otimes K^{n-r}
  & \longrightarrow &\L_\circ^{n-r}E \\
  & e_1 \otimes e_2 \otimes \omega & \longmapsto & 
  e_1 \otimes \widetilde pr_K(e_2 \otimes \omega) & \longmapsto &
  \s_E(e_1,e_2)\omega -\frac{1}{n-r+1}e_1^\# \es e_2 \dn \omega \\
  &&&&& \quad +\frac{r+2}{(n+r+3)(r+1)}e_1 \dn e_2^\# \es \omega \, .
\end{array}
$$
On the left side of (\ref{gl2}), $E \otimes E$ is seen as space of
endomorphisms, acting as derivations on $\L_\circ^{n-r}E$:
$$
  E \otimes E \cong \Sym^2E \oplus \L_\circ^2E \oplus \C \, ,
$$
where $\L_\circ^2E$ resp.~$\C$ denotes the trace--free part resp.~the
trace part of $\L^2E$. To calculate the trace--free part of an element
$ e_1 \wedge e_2 \in \L^2E$, one observes that the action on $E$ is
given by
$$
  (e_1 \wedge e_2)(e) = \s_E(e_1,e)e_2 - \s_E(e_2,e)e_1 \, ,
$$
its trace being $2 \s_E(e_1,e_2)$. Its
trace--free part is thus $e_1 \wedge e_2 -\frac{1}{n}\s_E(e_1,e_2)\id$.
Extended
as derivation to $\L_\circ^{n-r}E$, the identity $\id$ acts as multiplication
with $n-r$, and {\it a forteriori}, the trace--free part of $e_1 \wedge e_2$ 
becomes the operator 
$e_2\dn e_1^\# \es - e_1\dn e_2^\# \es - \frac{n-r}{n}\s_E(e_1,e_2)$.

With the preceeding calculations, the form of the projections is obvious:
$$
\begin{array}{rrcl}
  pr_{\C}:& E \otimes E \otimes \L_\circ^{n-r}E & 
\longrightarrow & \L_\circ^{n-r}E \\
   &e_1 \otimes e_2 \otimes \omega & \longmapsto & \s_E(e_1,e_2) \omega \, ,\\
\\
  pr_{\Sym^2E}:& E \otimes E \otimes \L_\circ^{n-r}E & 
\longrightarrow & \L_\circ^{n-r}E \\
   & e_1 \otimes e_2 \otimes \omega & \longmapsto & 
   e_2\dn e_1^\# \es\omega + e_1\dn e_2^\# \es\omega \\\\
   &&& = \big( \frac{r+2}{r+1} e_1\dn e_2^\# \es - e_1^\# \es e_2\dn +
     \s_E(e_1,e_2)\big)\omega \, ,\\
\\
  pr_{\L_\circ^2E}:& E \otimes E \otimes \L_\circ^{n-r}E & 
\longrightarrow & \L_\circ^{n-r}E \\
   & e_1 \otimes e_2 \otimes \omega & \longmapsto & 
   e_2\dn e_1^\# \es\omega - e_1\dn e_2^\# \es\omega - 
\frac{n-r}{n}\s_E(e_1,e_2)\omega \\\\
   &&& = \big(-\frac{r}{r+1}e_1\dn e_2^\# \es - e_1^\# \es e_2 \dn
    + \frac{r}{n}\s_E(e_1,e_2)\big)\omega.
\end{array}
$$
To summarize all calculations above let's consider a diagram:
\begin{center}
\begin{picture}(400,170)(-155,-60)

\put(-140,100){$\lb\C\otimes\L_\circ^{n-r}E\rb_{\L_\circ^{n-r}E}$}
\put(-140,85){$\s_E(e_1,e_2) \omega$}
\put(135,100){$\L_\circ^{n-r}E$}
\put(135,85){$e_1^\#\es e_2\dn$}
\put(135,40){$\L_\circ^{n-r}E$}
\put(135,25){$e_1\dn e_2^\#\es \omega$}
\put(0,40){$E \otimes E \otimes \L_\circ^{n-r}E$}
\put(0,25){$e_1 \otimes e_2 \otimes \omega$}
\put(135,-20){$\lb E \otimes K^{n-r}\rb_{\L_\circ^{n-r}E}$}
\put(110,-38){$\s_E(e_1,e_2)\omega -\frac{1}{n-r+1}e_1^\# \es e_2 \dn \omega$}
\put(110,-56){$+\frac{r+2}{(n+r+3)(r+1)}e_1 \dn e_2^\# \es \omega$}
\put(-165,-20){$\lb\L_\circ^2E\otimes\L_\circ^{n-r}E\rb_{\L_\circ^{n-r}E}$}
\put(-165,-38){$e_2\dn e_1^\#\es\omega - e_1\dn e_2^\#\es\omega$}
\put(-165,-56){$-\frac{n-r}{n}\s_E(e_1,e_2)\omega$}
\put(-165,40){$\lb\Sym^2E\otimes\L_\circ^{n-r}E\rb_{\L_\circ^{n-r}E}$}
\put(-165,22){$(e_2\dn e_1^\# \es + e_1\dn e_2^\#\es)\omega$}
\put(-10,45){\vector(-1,1){40}}
\put(-10,37){\vector(-1,0){40}}
\put(-10,30){\vector(-1,-1){40}}
\put(80,45){\vector(1,1){40}}
\put(80,37){\vector(1,0){40}}
\put(80,30){\vector(1,-1){40}}
\put(-25,70){$pr_{\C}$}
\put(-50,45){$pr_{\Sym^2E}$}
\put(-55,10){$pr_{\L^2_\circ}$}
\put(80,70){$pr_{-+}$}
\put(95,45){$pr_{+-}$}
\put(105,10){$pr_{K}$}
\end{picture}
\end{center}
Finally, we are able to give the $E$--part of the Weitzenb\"ock 
matrix ${\cal W}_E$:
\be
  {\cal W}_E = 
\left(\begin{array}{ccc}
  \frac{1}{n-r+1}           & -\frac{r+2}{(n+r+3)(r+1)}           & 1 \\\\
  -\frac{n-r}{n-r+1}        &  \frac{(n+r+2)(r+2)}{(n+r+3)(r+1)}  & 1\\\\
-\frac{(n-r)(n+1)}{n(n-r+1)}& -\frac{r(n+r+2)(n+1)}{n(n+r+3)(r+1)}& \frac{r}{n}
\end{array}\right) \, .
\ee

The full Weitzenb\"ock matrix is a $ 6\times 6$--matrix 
which is calculated as
Kronecker product $ {\cal W}={\cal W}_E \otimes {\cal W}_H $ of the
two partial Weitzenb\"ock matrices, i.~e.~the nine $2\times 2$-blocks
of ${\cal W}$ are obtained by multiplying $ {\cal W}_H$ with the corresponding
entry in ${\cal W}_E$. 
\be
\begin{array}{c}
{\cal W} = 
\left(\begin{array}{cccccc}
  \frac{1}{n-r+1} & -\frac{r}{(n-r+1)(r+1)} & -\frac{r+2}{(n+r+3)(r+1)} 
  & \frac{r(r+2)}{(n+r+3)(r+1)^2}  & 1&-\frac{r}{r+1} \\\\
  \frac{r}{n-r+1} & \frac{r(r+2)}{(n-r+1)(r+1)} 
  & -\frac{r(r+2)}{(n+r+3)(r+1)} & -\frac{r(r+2)^2}{(n+r+3)(r+1)^2}
  & r & \frac{r(r+2)}{r+1} \\\\
  -\frac{n-r}{n-r+1} & \frac{r(n-r)}{(n-r+1)(r+1)} 
  & \frac{(n+r+2)(r+2)}{(n+r+3)(r+1)} & -\frac{r(n+r+2)(r+2)}{(n+r+3)(r+1)^2}
  & 1 & -\frac{r}{r+1} \\\\
  -\frac{(n-r)r}{n-r+1} & -\frac{r(r+2)(n-r)}{(n-r+1)(r+1)} 
  & \frac{r(n+r+2)(r+2)}{(n+r+3)(r+1)} & \frac{r(n+r+2)(r+2)^2}{(n+r+3)(r+1)^2}
  & r & \frac{r(r+2)}{r+1} \\\\
  -\frac{(n-r)(n+1)}{n(n-r+1)} & \frac{r(n-r)(n+1)}{n(n-r+1)(r+1)} &
 -\frac{r(n+r+2)(n+1)}{n(n+r+3)(r+1)}  & 
   \frac{r^2(n+r+2)(n+1)}{n(n+r+3)(r+1)^2} 
  & \frac{r}{n} & -\frac{r^2}{n(r+1)} \\\\
  -\frac{r(n-r)(n+1)}{n(n-r+1)} & -\frac{r(r+2)(n-r)(n+1)}{n(n-r+1)(r+1)} 
&
-\frac{r^2(n+r+2)(n+1)}{n(n+r+3)(r+1)}  
& -\frac{r^2(r+2)(n+r+2)(n+1)}{n(n+r+3)(r+1)^2} 
  & \frac{r^2}{n} & \frac{r^2(r+2)}{n(r+1)}
\end{array}\right)\\ \hfil
\end{array}
\ee

\subsection{Associated Differential Operators}

In this section the two different sets of projectors onto copies of
$\spb(M)$ related by the associativity isomorphism 
\be \label{assi}
  (TM \otimes TM) \otimes \spb(M) \cong TM \otimes (TM \otimes \spb(M))
\ee
are applied to the special section 
$\nb^2 \psi\in \G(TM \otimes TM \otimes \spb(M))$.
This defines two sets of differential operators on sections of the spinor
bundle related by the matrix ${\cal W}$. In the sequel we determine all
these differential operators.

On the right hand side of (\ref{assi})
four of the projections of $\nb^2\psi$ are easily identified
with expressions like $\frac{1}{2}D^-_+D^+_-\psi$.
Due to the definition of the Clifford multiplication, the
explicit expressions have an additional factor $\frac{1}{2}$.
$$
\begin{array}{rrcl}
pr_{-+}\otimes pr_{-+}: &
  \nb^2\psi & \longmapsto & \sum_{i,j=0}^{2n}\sum_{a,b=1}^2
  (dh_a\esn\otimes de_i \es)(dh_b^\b\cdot\otimes de_j^\b \dn)
   \nb^2_{h_a\otimes e_i,h_b\otimes e_j}\psi \\[1ex]
  &&& = \frac{1}{2} D_-^-D_+^+\psi
      = -\frac{1}{2} (D_+^+)^\ast D_+^+\psi \, ,\\\\
pr_{+-}\otimes pr_{-+}: &
  \nb^2\psi & \longmapsto & \sum_{i,j=0}^{2n}\sum_{a,b=1}^2
  (dh_a^\b\cdot\otimes de_i \es)(dh_b\esn\otimes de_j^\b \dn)
   \nb^2_{h_a\otimes e_i,h_b\otimes e_j}\psi \\[1ex]
  &&& = \frac{1}{2} D_-^+D_+^-\psi \, ,\\\\ 
pr_{-+}\otimes pr_{+-}: &
  \nb^2\psi & \longmapsto & \sum_{i,j=0}^{2n}\sum_{a,b=1}^2
  (dh_a\esn\otimes de_i^\b \dn)(dh_b^\b\cdot\otimes de_j \es)
   \nb^2_{h_a\otimes e_i,h_b\otimes e_j}\psi \\[1ex]
  &&& = \frac{1}{2} D_+^-D_-^+\psi \, ,\\\\ 
pr_{+-}\otimes pr_{+-}: &
  \nb^2\psi & \longmapsto & \sum_{i,j=0}^{2n}\sum_{a,b=1}^2
  (dh_a^\b\cdot\otimes de_i^\b \dn)(dh_b\esn\otimes de_j \es)
   \nb^2_{h_a\otimes e_i,h_b\otimes e_j}\psi \\[1ex]
  &&& = \frac{1}{2} D_+^+D_-^-\psi
      = -\frac{1}{2} (D_-^-)^\ast D_-^-\psi \, .
\end{array}
$$
This is only a particular case of a more general phenomenon:
\begin{Lemma}
Let $V$ and $\S$ be the representations spaces for $TM$ and $\spb(M)$ and
let $p: V \otimes \S \to W$ be an $\Sp(n)\Sp(1)$--equivariant mapping.
Consider the equivariant mapping $ \Phi:  V \otimes (V \otimes \S) \to \S $
given by $ \id \otimes p^\ast p$ followed by contraction of the two
$V$--factors. Then the differential operator $\Phi \circ \nb^2 $
is equal to $ - P^\ast P$, where $P = p \circ \nb$.
\end{Lemma}

This lemma can be applied to the remaining two projectors.
To begin with $\Phi = pr_{-+}\otimes pr_{K}$ we consider the multiplication
$ m: H \otimes \Sym^rH \to \Sym^{r+1}H $. Its adjoint $m^\ast$ is given by
$ \sum h_i \otimes dh_i \esn$. The adjoint of the projection
$ \widetilde{pr}_K : E \otimes \L_\circ^{n-r}E \to K^{n-r}$ is simply the
inclusion. By definition, 
$pr_K: E \otimes E \otimes \L_\circ^{n-r}E \to \L_\circ^{n-r}E$ is 
$ \id_E \otimes \widetilde{pr}_K $ followed by contraction of the two 
$E$-factors. Similarly, it is easy to show that
$pr_{-+}: H \otimes H \otimes \Sym^rH \to \Sym^rH$ is $\id_H \otimes m^\ast m $
followed by contraction. The assumption of the lemma is thus satisfied and we
conclude that $ (pr_{-+}\otimes pr_K) \circ \nb^2 = -(T^+)^\ast T^+$.

The same argument goes through for the projector $pr_{+-}\otimes pr_{K}$. The 
adjoint of the contraction operator $ c: H \otimes \Sym^rH \to \Sym^{r-1}H $,
however, is given by $c^\ast = -\sum h_i \otimes dh_i^\b \cdot $. Hence,
$ (pr_{+-}\otimes pr_K) \circ \nb^2 = (T^-)^\ast T^-$.

The computed projections of $\nb^2\psi$ corresponding to the splitting on the
right hand side of (\ref{assi}) can be collected as entries of a vector:
\be
\Big(-\tfrac{1}{2}(D^+_+)^\ast D^+_+\psi \;\,~ \tfrac{1}{2}D^+_-D^-_+\psi \;\,~
    \tfrac{1}{2}D^-_+D^+_-\psi \;\,~ -\tfrac{1}{2}(D^-_-)^\ast D^-_-\psi\;\,~
   -(T^+)^\ast T^+\psi \;\,~ (T^-)^\ast T^-\psi \Big) \, .
\ee

To determine the differential operators occuring on the left hand 
side of the isomorphism, some technical lemmata are needed.
\begin{Lemma}
$$
\sum_{i,j=0}^{2n}\sum_{a,b=1}^2
 \big((dh_b^\b\cdot dh_a\esn + dh_a^\b\cdot dh_b\esn) \otimes
   \s_E(de_i^\b,de_j^\b)\big)
 \nb^2_{h_a\otimes e_i,h_b\otimes e_j}\psi
= \tfrac{r(r+2)}{n+2}\tfrac{\k}{4}\psi \, .
$$
\end{Lemma}
\leer
\proof
The left hand side factorizes over the projection
of $\nb^2\psi \in \G(TM\otimes TM \otimes\spb(M))$ onto
$\G(\Sym^2H \otimes \C \otimes \spb(M) \subset
\L^2TM \otimes \spb(M))$. Hence it is a curvature term,
and only the $H$--symmetric part $\tfrac{-\k}{8n(n+2)} R_H$
contributes according to Lemma \ref{curvature}.
$$
\begin{array}{rcl}
 \lefteqn{\sum_{i,j=0}^{2n}\sum_{a,b=1}^2 
 \big((dh_b^\b\cdot dh_a\esn + dh_a^\b\cdot dh_b\esn) \otimes
   \s_E(de_i^\b,de_j^\b)\big)
 \nb^2_{h_a\otimes e_i,h_b\otimes e_j}} &&\\\\
&=& \tfrac{1}{2}\displaystyle\sum_{i,j=0}^{2n}\sum_{a,b=1}^2 
\big((dh_b^\b\cdot dh_a\esn + dh_a^\b\cdot dh_b\esn)
\otimes \s_E(de_i^\b,de_j^\b)\big)
  \\\\
&&  \cdot
  \tfrac{-\k}{8n(n+2)}\big((h_b\cdot h_a^\#\esn + h_a\cdot h_b^\#\esn)
\otimes \s_E(e_i,e_j)\big)
   \\\\
&=& \tfrac{-\k}{4(n+2)}\displaystyle\sum_{a,b=1}^2
    dh_a^\b\cdot dh_b\esn(h_a\cdot h_b^\#\esn + h_b\cdot h_a^\#\esn) \\\\
&=& \tfrac{-\k}{4(n+2)}\displaystyle\sum_{a,b=1}^2
    \big(\s_H(dh_b^\b,h_a)dh_a^\b\cdot h_b^\#\esn 
+ dh_a^\b\cdot h_a\cdot dh_b\esn h_b^\#\esn
    \\\\
&& + dh_a^\b\cdot h_b\cdot dh_b\esn h_a^\#\esn 
+\s_H(dh_b^\b,h_b)dh_a^\b\cdot h_a^\#\esn\big)
    \\\\
&=& \tfrac{-\k}{4(n+2)}\big(-r -(r-1)r -2r\big) \\\\
&=& \tfrac{r(r+2)}{n+2}\tfrac{\k}{4} \, . \qed
\end{array}
$$

\begin{Lemma}
$$
\sum_{i,j=0}^{2n}\sum_{a,b=1}^2
 \big(\s_H(dh_a^\b,dh_b^\b) \otimes
   (de_j^\b \dn de_i\es + de_i^\b \dn de_j\es)\big)
 \nb^2_{h_a\otimes e_i,h_b\otimes e_j}
= \tfrac{(n+r+2)(n-r)}{n(n+2)}\tfrac{\k}{4} \, .
$$
\end{Lemma}
\leer
\proof
The situation here is exactly the same as in the preceeding lemma
only with the roles of $H$ and $E$ interchanged. The result only
depends on the $E$--symmetric part $\tfrac{-\k}{8n(n+2)} R_E$ because
$R^{hyper}$ does not contribute as has been shown in Corollary 
\ref{qzero}.
$$
\begin{array}{rcl}
  \lefteqn{\sum_{i,j=0}^{2n}\sum_{a,b=1}^2
   \big(\s_H(dh_a^\b,dh_b^\b) \otimes
   (de_j^\b \dn de_i\es + de_i^\b \dn de_j\es)\big)
   \nb^2_{h_a\otimes e_i,h_b\otimes e_j}} && \\\\
&=& \tfrac{1}{2}\sum_{i,j=0}^{2n}\sum_{a,b=1}^2
   \big(\s_H(dh_a^\b,dh_b^\b) \otimes (de_j^\b \dn de_i\es 
+ de_i^\b \dn de_j\es)\big)
   \\\\
&& \cdot
   \tfrac{-\k}{8n(n+2)}0
   \big(\s_H(h_a,h_b) \otimes (e_j \dn e_i^\#\es + e_i \dn e_j^\#\es)\big)
   \\\\
&=& \tfrac{-\k}{4n(n+2)}\sum_{i,j=0}^{2n}
    de_i^\b \dn de_j\es(e_i \dn e_j^\#\es + e_j \dn e_i^\#\es)\\\\
&=& \tfrac{-\k}{4n(n+2)}\sum_{i,j=0}^{2n}
    \Big(de_i^\b \dn e_i \dn de_j\es e_j^\# \es 
+ \s_E(de_j^\b,e_i)de_i^\b \dn e_j^\# \es
    + \tfrac{1}{r+1} de_i^\b \dn de_j\dn e_i^\# \es e_j^\#\es \\\\
&&   - de_i^\b \dn e_j \dn de_j\es e_i^\# \es 
+ \s_E(de_j^\b,e_j)de_i^\b \dn e_i^\# \es
    + \tfrac{1}{r+1} de_i^\b \dn de_j\dn e_j^\# \es e_i^\#\es \Big) 
     \\\\
&=& \tfrac{-\k}{4n(n+2)}\sum_{i,j=0}^{2n}
    \big(-(n-r) + (n-r-1)(n-r) - 2n(n-r)\big) \\\\
&=& \tfrac{(n+r+2)(n-r)}{n(n+2)}\tfrac{\k}{4} \, . \qed
\end{array}
$$

\begin{Lemma}
$$
  pr_{\Sym^2H}\otimes pr_{\L_\circ^2E}(\nb^2 \psi) = 0 \, .
$$
\end{Lemma}
\leer
\proof
The left hand side is the projection of $\nb^2\psi$ onto
$ \G(\Sym^2H \otimes \L_\circ^2E \otimes \spb(M)) \subset 
\G(\L^2TM \otimes \spb(M))$. Hence, it is again a part of the
curvature tensor of the manifold. But Lemma \ref{curvature} shows
that this contribution does not exist.
\leer
The rest is easy. With help of the preceeding lemmata, 
the projections of $\nb^2\psi$ onto the irreducible $\spb(M)$--summands of
the left hand side of $(TM \otimes TM)
\otimes \spb(M) \cong TM \otimes (TM \otimes\spb(M))$
can be given:
$$
\begin{array}{rrcl}
pr_{\C}\otimes pr_{\C}: &
  \nb^2\psi & \longmapsto & \sum_{i,j=0}^{2n}\sum_{a,b=1}^2
  \s_H(dh_a^\b,dh_b^\b)\s_E(de_i^\b,de_j^\b)
   \nb^2_{h_a\otimes e_i,h_b\otimes e_j}\psi \\\\
  &&& = -\nb^\ast\nb\psi \, ,\\\\ 
pr_{\Sym^2H}\otimes pr_{\C}: &
  \nb^2\psi & \longmapsto & 
   \frac{r(r+2)}{n+2}\frac{\k}{4} \psi \, ,\\\\ 
pr_{\C}\otimes pr_{\Sym^2E}: &
  \nb^2\psi & \longmapsto & 
  \frac{(n+r+2)(n-r)}{n(n+2)}\frac{\k}{4} \psi\, ,\\\\ 
pr_{\Sym^2H}\otimes pr_{\Sym^2E}: &
  \nb^2\psi & \longmapsto & \sum_{i,j=0}^{2n}\sum_{a,b=1}^2
  (dh_b^\b\cdot dh_a\esn + dh_a^\b\cdot dh_b\esn) \\\\
  &&&\otimes
  (de_j^\b \dn de_i\es + de_i^\b \dn de_j\es)
   \nb^2_{h_a\otimes e_i,h_b\otimes e_j}\psi \\\\
  &&& =: {\cal C}\psi \, , \\\\
pr_{\C}\otimes pr_{\L_\circ^2E}: &
  \nb^2\psi & \longmapsto & \sum_{i,j=0}^{2n}\sum_{a,b=1}^2
  \s_H(dh_a^\b,dh_b^\b) \\\\
  &&&\otimes 
  \big(de_j^\b\dn de_i \es\omega - de_i^\b\dn de_j \es\omega 
    - \frac{n-r}{n}\s_E(de_i^\b,de_j^\b) \big)\psi \\\\
  &&& =: {\cal L}\psi \, ,\\\\
pr_{\Sym^2H}\otimes pr_{\L_\circ^2E}: &
  \nb^2\psi & \longmapsto & 0 \, .
\end{array}
$$
Here, ${\cal C}$ and ${\cal L}$ are mysterious 2nd order differential
operators which are simply defined by the expressions above.
As for the right hand the computed projections can be collected in
a vector:
\be
  \Big( -\nb^\ast\nb\psi\;\,~ \tfrac{r(r+2)}{n+2}\tfrac{\k}{4}\psi\;\,~
  \tfrac{(n+r+2)(n-r)}{n(n+2)}\tfrac{\k}{4}\psi\;\, ~{\cal C}\psi\;\,~
  {\cal L}\psi\;\,~0 \Big) \, .
\ee

Summarizing all results, the full Weitzenb\"ock matrix equation reads
as follows:
\be \label{last}
\left(
\begin{array}{c}
  -\nb^\ast\nb\psi \\\\ \frac{r(r+2)}{n+2}\frac{\k}{4}\psi \\\\
  \frac{(n+r+2)(n-r)}{n(n+2)}\frac{\k}{4}\psi  \\\\ {\cal C}\psi \\\\
  {\cal L}\psi \\\\  0
\end{array}
\right)
= {\cal W} \cdot
\left(
\begin{array}{c}
  -\frac{1}{2}(D^+_+)^\ast D^+_+\psi \\\\ \frac{1}{2}D^+_-D^-_+\psi \\\\
    \frac{1}{2}D^-_+D^+_-\psi \\\\ -\frac{1}{2}(D^-_-)^\ast D^-_-\psi \\\\
   -(T^+)^\ast T^+\psi \\\\ (T^-)^\ast T^-\psi
\end{array}
\right) \, .
\ee

\section{Proof of the Theorem}

The proof of the lower bound for the Dirac spectrum
is an easy application of our
Weitzenb\"ock matrix $ {\cal W} $. Nevertheless, before deriving this result
we want to give a few other consequences showing how much
information is contained in  $ {\cal W} $.

From equation \eqref{last} it is clear that the matrix  $ {\cal W} $
yields six lineary independent Weitzenb\"ock  formulas involving
Dirac and twistor
operators. Since we have no further information about the operators
$ {\cal C } $   and  $ {\cal L } $, we can use only four of them and
look for suitable linear combinations to eliminate certain operators.
This is easily done by multiplying equation
\eqref{last} from the left by a row vector
$ (a_1 \; \ldots\; a_6) $, where the $ a_j $ are arbitrary
coefficients. For example, we can eliminate the two twistor operators
$ T^{\pm}  $ by multiplying the second equation with
$ \; \tfrac{r}{n} \; $ and subtracting the last. This is obtained by
multiplying with
$ (0\;\,~ \tfrac{r}{n}\;\,~ 0\;\,~ 0\;\,~ 0\;\,~ -1) $,
and thus the corresponding Weitzenb\"ock formula can be written as
\be\label{nn}
\tfrac{r}{2} \,
D^-_- \,  D^+_+
\; + \;
\tfrac{r(r+2)}{2(r+1)} \,
D^+_- \,  D^-_+
\; + \;
\tfrac{r^2}{2(r+1)} \,
D^-_+ \,  D^+_-
\;  + \;
\tfrac{r^2(r+2)}{2(r+1)^2} \,
D^+_+ \, D^-_-
\; = \;
\tfrac{r^2(r+2)}{n(n+2)} \,  \tfrac{\kappa}{4} \,         .
\ee
This equation also follows from a theorem of Y.~Nagatomo and T.~Nitta.
In \cite{naga} they consider Dirac and twistor operators
on arbitrary bundles   $ S^p H \otimes \Lambda^q E $  and prove
a general Weitzenb\"ock formula.

Since we have four independent equations for six unknowns
($ D^-_-  D^+_+  \psi$, $D^+_-   D^-_+  \psi$ etc.)
it is possible to eliminate all four twistor operators
at once. The corresponding row vector for $r \neq 0$ is
$ ( -1\;\,~ \tfrac{1}{n}\;\,~ 1\;\,~ 0\;\,~ 0\;\,~ -\tfrac{1}{r}) $,
and we obtain the well--known Lichnerowicz formula:
$$
 \, \nabla^* \nabla    \, + \, \tfrac{\kappa}{4} \,
\; = \;
\, D^+_- \,  D^-_+
\, + \,
D^-_+ \,  D^+_-
\; = \;
D^2   \, .
$$

Before coming to the most important application of this 
formalism, the proof of the eigenvalue estimate, let's have
a look on the particular form of eigenspinors of $D$.
Since $D^2$ respects the splitting into the subbundles $\spb_r(M)$,
some assumptions on the form of an eigenspinor of $D$ can be made.
Let $\psi_r \in \G(\spb_r(M))$ be an eigenspinor for $ D^2 $
with eigenvalue $ \lambda^2 $. By $D$, it is mapped to
$ D\psi_r =: \psi_{r-1} + \psi_{r+1} \in \G(\spb_{r-1}(M)\oplus
\spb_{r+1}(M))$. Because of $D^2 \psi_r \in \G(\spb_r(M))$ we get
$D^-_+ \psi_{r-1} = 0 = D^+_- \psi_{r+1}$. But starting
with $\psi_{r-1}$, which itself is an eigenspinor of $D^2$ with the
eigenvalue $\l^2$, we see that $\l \psi_{r-1} \pm D\psi_{r-1}
\in \G(\spb_{r-1}(M)\oplus \spb_r(M))$ is an eigenspinor of $D$ with
eigenvalue $\pm \l$. The analogous argument goes through with $\psi_{r+1}$.
Hence we always can assume that an eigenspinor for $D$ is of the form
$\psi = \psi_r + \psi_{r+1} \in \G(\spb_r(M)\oplus \spb_{r+1}(M))$.

\begin{Theorem}
Let $( M^{4n}, \, g) $ be a compact quaternionic K\"ahler spin
manifold of positive scalar curvature $ \kappa $ and let
$\psi = \psi_r + \psi_{r+1}\in \G(\spb_r(M)\oplus \spb_{r+1}(M))$
be an eigenspinor for $D$ with eigenvalue $\l$. Then
$$
\lambda^2
\; \ge  \;
\frac{n + r + 3}{n + 2 } \,  \frac{\kappa}{4} \, .
$$
\end{Theorem}
\leer
\proof
We apply formula (\ref{last}) to the spinor $\psi_r$ and multiply
it from the left with the row vector
$ (0\;\,~ \tfrac{n+r+2}{n}\;\,~ r+2\;\,~ 0\;\,~ 0\;\,~ -\tfrac{r+2}{r}) $.
This eliminates the operators $T^-$ and $D^-_-$ and we
find after taking the $L^2$-product with $\psi_r$.
$$
\begin{array}{c}
- \, \tfrac{r+1}{n-r+1} \,  \| D^+_+ \psi_r \|^2
\; + \;
(r+2) \,  \| D^-_+ \psi_r \|^2
\; + \;
\tfrac{(r+2)(n+r+2)}{n+r+3} \,  \| D^+_- \psi_r \|^2
\; - \;
2(r+1) \,  \| T^+ \psi_r \|^2 \\\\
 = \;
\tfrac{(r+2)(n+r+2)}{n+2} \, \tfrac{\kappa}{4} \,  \| \psi_r \|^2 \,.
\end{array}
$$
Using $D^-_+ \psi_r = 0$ and
estimating   $  \; \| D^+_+ \psi \|^2  \; $  and
$ \;    \| T^+ \psi \|^2   \; $ by zero leads to the inequality
$$
\| D^+_- \psi_r \|^2
\; \ge  \;
\tfrac{n+r+3}{n+2} \, \tfrac{\kappa}{4} \,  \| \psi_r \|^2  \, .
$$
The same procedure is carried out with $\psi_{r+1}$. The only difference
is that in this case $D^+_-\psi_{r+1}=0$, hence the same inequality holds:
$$
\| D^+_- \psi_{r+1} \|^2
\; \ge  \;
\tfrac{n+r+3}{n+2} \, \tfrac{\kappa}{4} \,  \| \psi_{r+1} \|^2  \, .\qed
$$


\begin{thebibliography}{MMM99-9}

\bibitem[Ale68-1]{alex2}
  Alekseevskii, D.V.,
  Riemannian spaces with exceptional holonomy groups,
  Funktional Anal. Appl. 2, 1968, 97--105;

\bibitem[Ale68-2]{alex2a}
  Alekseevskii, D.V.,
  Compact quaternion spaces,
  Funktional Anal. Appl. 2, 1968, 106--114;

\bibitem[Bon67]{bon}
  Bonan, E.,
  Sur les {G}--structures de type quaternionien,
  Cahiers Top. et Geom. Diff. 9, 1967, 389--461;

\bibitem[Fri80]{fri5}
  Friedrich, Th.,
  Der erste {E}igenwert des {D}irac--{O}perators einer kompakten
  {R}iemannschen {M}annigfaltigkeit nichtnegativer {S}kalarkr\"ummung,
  Math. Nachr. 97, 1980, 117--146;

\bibitem[Hij86]{hij2}
  Hijazi, O.,
  A conformal lower bound for the smallest
  eigenvalue of the {D}irac operator and {K}illing spinors,
  Comm. Math. Phys. 104, 151--162;

\bibitem[Hij94]{hij1}
  Hijazi, O.,
  Eigenvalues of the {D}irac operator on {C}ompact
  {K}\"ahler {M}anifolds,
  Comm. Math. Phys. 160, 1994, 563--579;

\bibitem[Hij96]{hij}
  Hijazi, O.,
  Twistor operators and eigenvalues of the
              {D}irac operator on compact quaternionic--{K}\"ahler
               spin manifolds, 
Preprint,  1996;

\bibitem[HiM95]{hij3}
  Hijazi, O., Milhorat, J.--L.,
  Minoration des valeurs propres de
  de l'op\'{e}rateur de {D}irac sur les vari\'{e}t\'{e}s
  spin {K}\"ahler--quaternioniennes,
  J. Math. Pures Appl. 74, 1995, 387--414;

\bibitem[Ish74]{ish3}
  Ishihara, S.,
  Quaternionic {K}\"ahlerian manifolds,
  J. Diff. Geom. 9, 1974, 483--500;

\bibitem[Kir86]{kir2}
  Kirchberg, K.-D.,
  An estimation for the first eigenvalue of the {D}irac operator on
  closed {K}\"ahler manifolds of positive scalar curvature,
  Ann. Global Anal. Geom. 4, 1986, 291--325;

\bibitem[Kir90]{kir6}
  Kirchberg, K.-D.,
  The first eigenvalue of the {D}irac operator   on
                 {K}\"ahler manifolds,
  J. Geom. Phys. 7, 1990, 449--468;

\bibitem[Kra66]{krain}
  Kraines, V.Y.,
   Topology of quaternionic manifolds,
  Trans. Amer. Math. Soc. 122, 1966, 357--367;

\bibitem[Lic87]{lic}
   Lichnerowicz, A.,
   Spin manifolds, {K}illing spinors and the universality of the
   Hijazi inequality,
   Letters in Math. Phys. 13, 1987, 331--344;
 
\bibitem[Lic90]{lic1}
   Lichnerowicz, A.,
  La premi\`{e}re valeur propre de l'op\'{e}rateur de {D}irac
  pour une vari\`{e}t\`{e}s {K}\"ahl\`{e}rienne et son cas limite,
   C.R.Acad. Sci. Paris 311, 1990, 717--722;

\bibitem[Mil92]{mil}
  Milhorat, J.--L.,  
  Spectre de l'op\'{e}rateur de {D}irac sur les espaces projectifs
  quaternioniens,
  C.R.Acad. Sci. Paris 314, 1992, 69--72;

\bibitem[Nag96]{naga}
  Nagatomo, Y. and Nitta, T.,
  Vanishing {T}heorems for {Q}uaternionic {C}omplexes,
  Preprint, 1996;

\bibitem[Sal82]{sala}
  Salamon, S.M.,
  Quaternionic {K}\"ahler manifolds,
  Invent. Math. 67, 1982, 143--171.

\end{thebibliography}
\end{document}